\documentclass[11pt]{amsart}   	
\usepackage{geometry}             
\usepackage[dvipsnames]{xcolor}   		             		
\usepackage{graphicx}				
\usepackage{mathrsfs}
\usepackage{amssymb,amsthm,amsmath,mathtools}
\usepackage{tikz}
\usepackage{tikz-cd}
\usepackage{accents,upgreek,enumerate}
\usepackage[headings]{fullpage}
\usepackage{bm}
\usepackage[all]{xy}
\usepackage{caption}
\usepackage{floatrow}       
\usepackage{csquotes}
\usepackage{enumitem}
\usepackage{comment}
\usepackage{rotating}       
\usepackage{makecell}       
\usepackage{multirow}       
\usepackage{crossreftools}  
\usepackage{mathdots}
\usepackage{hyphenat}

\usetikzlibrary{calc}
\usetikzlibrary{fadings}
\usetikzlibrary{decorations.pathmorphing}
\usetikzlibrary{decorations.pathreplacing}

\renewcommand{\familydefault}{ppl}
\DeclareMathAlphabet{\mathpzc}{OT1}{pzc}{m}{it}

\usepackage[backref]{hyperref}
\hypersetup{
colorlinks   = true,          
urlcolor     = violet,          
linkcolor    = Purple,          
citecolor   = RubineRed             
}
\usepackage[capitalize,noabbrev]{cleveref}

\newtheorem{theorem}{Theorem}[section]

\newtheorem{lemma}[theorem]{Lemma}
\newtheorem{proposition}[theorem]{Proposition}

\newtheorem{quasi-theorem}[theorem]{Quasi-Theorem}

\theoremstyle{definition}
\newtheorem{definition}[theorem]{Definition}

\newtheorem{example}[theorem]{Example}
\newtheorem{remark}[theorem]{Remark}

\newcommand{\direction}{\vec{\mathbf{v}}} 

\newcommand{\A}{{\mathbb{A}}}                     
            
\newcommand{\NN} {{\mathbb N}}		
         
\newcommand{\QQ} {{\mathbb Q}}	

\newcommand{\RR} {{\mathbb R}}	
\newcommand{\CC} {{\mathbb C}}
\newcommand{\ZZ} {{\mathbb Z}}	
\renewcommand{\AA}{\mathbb{A}}	

\def\setminus{\smallsetminus}

\newcommand{\Puiseux}[2][t]{#2\{\!\{#1\}\!\}}

\DeclareMathOperator{\image}{Im}

\DeclareMathOperator{\Tr}{Tr}               

\DeclareMathOperator{\chara}{char}

\DeclareMathOperator{\val}{val}

\DeclareMathOperator{\paths}{path}

\newcommand{\gw}[1]{\left\langle #1 \right\rangle}
\newcommand{\qinv}[1]{\left\langle#1\right\rangle}

\DeclareMathOperator{\GW}{GW}                   

\def\trop{\mathrm{trop}}

\def\mult{\mathrm{mult}}

\newcommand{\Trop}{\operatorname{Trop}}

\makeatletter
\def\blfootnote{\xdef\@thefnmark{}\@footnotetext}
\makeatother

\definecolor{darkspringgreen}{rgb}{0.09, 0.6, 0.1}
\definecolor{dsg}{rgb}{0.09, 0.6, 0.1}

\title[Tropical methods for counting plane curves]{\bf Tropical methods for counting plane curves --- complex, real and quadratically enriched} 

\date{}

\author {Andr\'es Jaramillo Puentes}
\address { Universit\`{a} di Catania, Dipartimento di Matematica e Informatica,
Viale Andrea Doria 6, 95125 Catania, Italy}
\email {ga.jaramillopuentes@unict.it}
\author {Hannah Markwig}
\address {Universit\"at T\"ubingen, Fachbereich Mathematik, Auf der Morgenstelle 10, 72076 T\"ubingen, Germany }
\email {hannah@math.uni-tuebingen.de}
\author {Sabrina Pauli}
\address {TU Darmstadt, Fachbereich Mathematik, Schlossgartenstra\ss{}e 7, 64289 Darmstadt, Germany}
\email {pauli@mathematik.tu-darmstadt.de}
\author{Felix R\"ohrle}
\address {Universit\"at T\"ubingen, Fachbereich Mathematik, Auf der Morgenstelle 10, 72076 T\"ubingen, Germany }
\email{roehrle@math.uni-tuebingen.de}

\subjclass[2020]{14N10, 14N35, 14T20, 14T25, 14P99}
\keywords{Gromov-Witten invariants, Welschinger invariants, tropical curves, quadratically enriched counts}

\begin{document}

\begin{abstract}
Since the first famous correspondence theorem by Mikhalkin appeared in 2005, tropical geometry has allowed a parallel treatment of real and complex counting problems. A prime example are the genus 0 Gromov-Witten invariants of the plane which count rational plane curves of degree~$d$ satisfying point conditions and their real counterpart, the Welschinger invariants, which both can be determined using tropical methods. Remarkably, the tropical computation of the two types of invariants works entirely in parallel. Recently, quadratically enriched enumerative geometry enables us to combine such real and complex counts under one roof, providing a simultaneous approach which can also be used for counts over other fields. Tropical geometry is a successful tool for the study and computation of such quadratically enriched enumerative invariants, too. In this survey, we provide an overview of tropical methods for plane curve counting problems over the real and complex numbers, and the new quadratically enriched counts.

\end{abstract}

\maketitle

\vspace{-0.2in}

\setcounter{tocdepth}{1}

\section{Introduction}
Tropical geometry can be viewed as a degeneration that replaces a plane algebraic curve by a piecewise linearly embedded graph in $\RR^2$ satisfying certain conditions. Even though the degeneration is severe, important geometric properties can often still be read off a tropicalized curve, thus enabling an exchange of information between algebraic geometry and combinatorics.

The foundational ideas for tropical geometry had appeared in different forms in the 1970s and 80s, e.g.\ as Viro's patchworking method \cite{Viro}, via amoebas \cite{EKL04}, Maslov's dequantization \cite{Lit05}, or logarithmic limit sets \cite{Ber71, BG84}. Only since the late 1990s an effort has been made to systematically exploit connections to questions in classical algebraic geometry, prominently represented in work of Mikhalkin and Sturmfels \cite{Mi03, Mi06, Stu02, SS04a}.

Even after it had become clear that tropicalization is a powerful tool for the study of algebraic geometry, the discovery that important enumerative invariants can be determined with tropical methods, came as a surprise. The first correspondence theorem, proved by Mikhalkin \cite{Mi03}, stating the equality of an enumerative invariant to its tropical counterpart, became a door opener for many applications of tropical methods in enumerative geometry. The tropical method has been particularly successful for problems in real enumerative geometry \cite{IKS03, IKS05, Mik15}, since, roughly speaking, the tropicalization does not depend on the underlying field (only the inverse operation, sometimes called realisation, i.e.\ the question which algebraic curves or how many tropicalize to a given tropical one, of course highly depends on the chosen ground field).

In recent years, two influential ways of combining real and complex enumerative counts under one roof appeared: on the one hand, the so-called refined invariants or Block-G\"ottsche invariants which discuss polynomial counts that specialize to the real or complex count after inserting particular values \cite{BG14}, on the other hand, the so-called quadratically enriched counts for which we equip our geometric objects to be counted with a quadratic form in such a way that the rank equals the complex count and the signature the real count \cite{KW17, LV19, Lev17, LevineWelschinger}. While both topics share this basic motivation, their general approaches are rather different. 

In this survey, based on the 2025 SRI talk by the second author, we focus on quadratically enriched counts of plane rational curves satisfying point conditions.
Our main motivation is to showcase tropical geometry as a tool for enumerative geometry, in the broad sense ranging from complex over real to quadratically enriched enumerative geometry.
The survey should serve as a gentle introduction to tropical curve counting methods and their applications. 

After we shortly present basics of tropical plane curves in Section \ref{sec-tropicalcurves} (including an overview of their relation to algebraic curves, i.e.\ a short introduction to tropicaliztion in Section \ref{sec-tropicalization}), we focus on the enumerative problem at the center of attention of this survey: the count of rational plane curves satisfying point conditions. 
We meet this enumerative problem in different guises: in Section \ref{sec-complex}, we discuss the complex version of the enumerative problem, yielding the so-called Gromov-Witten invariants of the plane. We state how the complex count can be expressed in terms of a tropical count in Section \ref{subsec-counttrop}. In tropical enumerative geometry, we do not just count tropical curves one by one, but we count with multiplicity which should reflect how many complex curves tropicalize to a given tropical curve. The way to determine this complex multiplicity is also discussed in Section \ref{subsec-counttrop}. In Section \ref{sec-corres}, we discuss the basic idea of the correspondence theorem for this complex count.

In Section \ref{sec-W}, we turn our attention to the real count, the Welschinger invariants. What is nice is that one can use exactly the same set of tropical curves to count, it is only the multiplicity that has to be adapted to the real situation, as discussed in Section \ref{sec-tropW}.

After introducing the background of quadratically enriched counts in Section \ref{sec-QEbasics}, in Section \ref{sec-QE}, we introduce the analogous quadratically enriched count of rational plane curves satisfying point conditions. Again, for the tropical analogue we can use the same set of tropical curves, only now with quadratically enriched multiplicities, as discussed in Section \ref{sec-tropQE}.
It is remarkable and beautiful that the same tropical curves (only with adapted multiplicities) can be used to obtain the complex, real and quadratically enriched count of rational plane curves satisfying $k$-rational point conditions (i.e.\ each point condition is defined over the base field).

Still, the question remains how one can actually determine a tropical count algorithmically, so that one can make proper use of the tropical method. This topic is discussed in Section \ref{sec-alg}. We focus in particular on a lattice path algorithm for counting tropical curves that has not yet appeared in the context of the quadratically enriched count. While it is obvious to experts how the complex and real versions of the algorithm have to be adapted to produce a quadratically enriched count, we view it as our small original contribution to this survey to actually spell out the details here. 

In Section \ref{sec-QEQ} finally, we discuss an exciting new extension of our enumerative problem: for the quadratically enriched count, we can also impose point conditions with points not defined over our ground field, but over a finite separable field extension. When specialized to the real version, this means counting real curves that do not only pass through real points but also through pairs of complex conjugate points. Tropical versions of this extended real count have been studied in \cite{Shu06b}. For arbitrary fields, the situation becomes more intricate, both in terms of the combinatorics underlying the tropical curves to be counted as well as in terms of the preimages under tropicalization that have to be studied. We give a quick overview over the recent methods to compute these numbers, and showcase an example with a tropical computation. In particular, we point out how many different counts these extended quadratically enriched enumerative invariants encompass: we can deduce as many as six special real and complex versions of this count from just one computation and even more quadratically enriched counts. Hoping that this example convinces the reader of the power of the tropical approach for quadratically enriched enumerative geometry, we end the survey with a short discussion of some open problems in the area which we feel deserve further attention in the future. 

\subsection{Acknowledgements}
The first, second and fourth author acknowledge support by DFG grant MA 4797/9-1. 
The third author acknowledges support by Deutsche Forschungsgemeinschaft (DFG, German Research Foundation) through the Colla\-borative Research Centre TRR 326 \textit{Geometry and Arithmetic of Uniformized Structures}, project number 444845124.

We all thank the SRI 2025 organizers for this special enormous conference event that provided so many occassions for stimulating discussions with colleagues.

\section{What is a tropical plane curve?}

There are many excellent introductions to tropical geometry which all discuss the basic case of tropical plane curves in detail (e.g.\ \cite{BIMS14, BS14, CMR23, IMS09, MS15, RST03}). To make this survey accessible to a broad audience, we introduce the main concepts shortly while trying to stay concise and refer to the literature for more details.

We first introduce tropical plane curves intrinsically, as certain graphs in the real plane (Section \ref{sec-tropicalcurves}), and then discuss in what sense they arise as degenerations of algebraic curves (Section \ref{sec-tropicalization}).

\subsection{Tropical plane curves as piecewise linear graphs dual to Newton subdivisions}\label{sec-tropicalcurves}

Roughly, a plane tropical curve is an edge weighted piecewise integer affine linearly embedded graph in $\mathbb{R}^2$ which satisfies the balancing condition at every vertex. The \emph{balancing condition} states that the weighted sum of the outward pointing primitive integer vectors in direction of the edges adjacent to a vertex must sum to $0$. The unbounded edges of a tropical curve are called \emph{ends}. 

\begin{example}

Figure \ref{fig-planecurves} shows two examples of tropical plane curves. The edges that are not adjacent to two vertices are supposed to be unbounded, i.e.\ ends. The $2$ in the left picture denotes the weight of the edge next to which it is written. All other weights are $1$ and left out in the picture. The balancing condition can be checked at every vertex. At the upper left vertex of the left picture for example, the three adjacent edges are of primitive direction $(-1,0)$, $(0,-1)$ and $(1,1)$ which adds up to $0$. At the left vertex of the edge of weight $2$, the three adjacent edges are of primitive direction $(-2,1)$, $(0,-1)$ and $(1,0)$. The weighted sum of these, $(-2,1)+(0,-1)+ 2 \cdot (1,0)$, is $0$.

    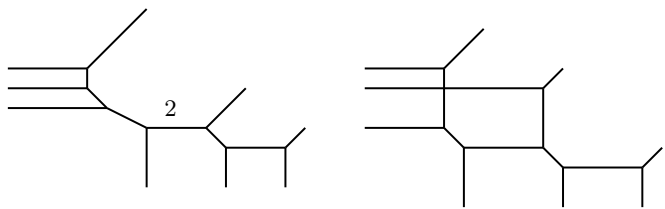
\begin{figure}
        \centering

\tikzset{every picture/.style={line width=0.75pt}} 

\begin{tikzpicture}[x=0.75pt,y=0.75pt,yscale=-1,xscale=1]

\draw    (100,70) -- (140,70) ;
\draw    (140,70) -- (170,40) ;
\draw    (140,80) -- (140,70) ;
\draw    (100,80) -- (140,80) ;
\draw    (140,80) -- (150,90) ;
\draw    (100,90) -- (150,90) ;
\draw    (150,90) -- (170,100) ;
\draw    (170,100) -- (200,100) ;
\draw    (170,130) -- (170,100) ;
\draw    (200,100) -- (220,80) ;
\draw    (200,100) -- (210,110) ;
\draw    (210,110) -- (240,110) ;
\draw    (210,130) -- (210,110) ;
\draw    (240,130) -- (240,110) ;
\draw    (240,110) -- (250,100) ;
\draw    (280,70) -- (320,70) ;
\draw    (320,100) -- (320,70) ;
\draw    (320,70) -- (340,50) ;
\draw    (280,100) -- (320,100) ;
\draw    (320,100) -- (330,110) ;
\draw    (330,110) -- (370,110) ;
\draw    (330,140) -- (330,110) ;
\draw    (370,80) -- (280,80) ;
\draw    (370,80) -- (380,70) ;
\draw    (370,110) -- (370,80) ;
\draw    (370,110) -- (380,120) ;
\draw    (380,140) -- (380,120) ;
\draw    (380,120) -- (420,120) ;
\draw    (420,140) -- (420,120) ;
\draw    (420,120) -- (430,110) ;

\draw (184.75,90.25) node  [font=\footnotesize] [align=left] {\begin{minipage}[lt]{8.67pt}\setlength\topsep{0pt}
$\displaystyle 2$
\end{minipage}};

\end{tikzpicture}

        \caption{Two tropical plane curves.}
        \label{fig-planecurves}
    \end{figure}
\end{example}

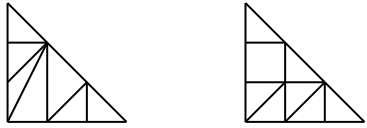
\begin{figure}
    \centering

\tikzset{every picture/.style={line width=0.75pt}} 

\begin{tikzpicture}[x=0.75pt,y=0.75pt,yscale=-1,xscale=1]

\draw    (350,120) -- (410,120) ;
\draw    (350,60) -- (410,120) ;
\draw    (370,120) -- (370,80) ;
\draw    (350,60) -- (350,120) ;
\draw    (390,120) -- (390,100) ;
\draw    (370,120) -- (390,100) ;
\draw    (370,100) -- (390,100) ;
\draw    (350,120) -- (370,100) ;
\draw    (350,100) -- (370,100) ;
\draw    (350,80) -- (370,80) ;
\draw    (230,120) -- (290,120) ;
\draw    (230,60) -- (290,120) ;
\draw    (250,120) -- (250,80) ;
\draw    (230,60) -- (230,120) ;
\draw    (270,120) -- (270,100) ;
\draw    (250,120) -- (270,100) ;
\draw    (230,100) -- (250,80) ;
\draw    (230,120) -- (250,80) ;
\draw    (230,80) -- (250,80) ;

\end{tikzpicture}

    \caption{The dual Newton subdivision of the two tropical plane curves depicted in Figure \ref{fig-planecurves}.}
    \label{fig-dual}
\end{figure}

The balancing condition implies that the weighted orthogonals of the outward pointing primitive vectors of the edges adjacent to a vertex (each turned in the same direction) form a closed path which coincides with the boundary edges of a polygon. Furthermore, the polygons of two adjacent vertices fit nicely together, and overall, we obtain a subdivision of a polygon dual to our plane tropical curve, called the \emph{dual Newton subdivision}.
Figure \ref{fig-dual} shows the dual Newton subdivisions for the two tropical plane curves depicted in Figure \ref{fig-planecurves}. 
Every vertex of the tropical plane curve is dual to a polygon in the subdivision, every edge to an orthogonal edge, every end to a boundary edge of the subdivision, and every connected component of $\mathbb{R}^2$ minus the tropical plane curve to a vertex of the dual subdivision. 
For more details on dual Newton subdivisions, see e.g. \cite{MS15, Mi03, RST03}.

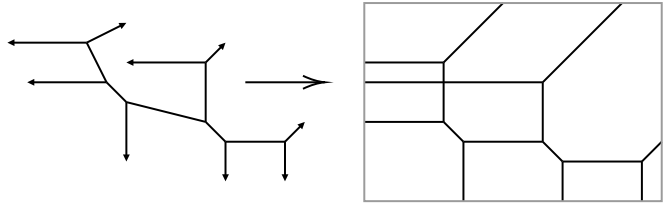
\begin{figure}
    \centering

\tikzset{every picture/.style={line width=0.75pt}} 

\begin{tikzpicture}[x=0.75pt,y=0.75pt,yscale=-1,xscale=1]

\draw    (163,70) -- (200,70) ;
\draw [shift={(160,70)}, rotate = 0] [fill={rgb, 255:red, 0; green, 0; blue, 0 }  ][line width=0.08]  [draw opacity=0] (3.57,-1.72) -- (0,0) -- (3.57,1.72) -- cycle    ;
\draw    (140,60) -- (157.32,51.34) ;
\draw [shift={(160,50)}, rotate = 153.43] [fill={rgb, 255:red, 0; green, 0; blue, 0 }  ][line width=0.08]  [draw opacity=0] (3.57,-1.72) -- (0,0) -- (3.57,1.72) -- cycle    ;
\draw    (103,60) -- (140,60) ;
\draw [shift={(100,60)}, rotate = 0] [fill={rgb, 255:red, 0; green, 0; blue, 0 }  ][line width=0.08]  [draw opacity=0] (3.57,-1.72) -- (0,0) -- (3.57,1.72) -- cycle    ;
\draw    (140,60) -- (150,80) ;
\draw    (113,80) -- (150,80) ;
\draw [shift={(110,80)}, rotate = 0] [fill={rgb, 255:red, 0; green, 0; blue, 0 }  ][line width=0.08]  [draw opacity=0] (3.57,-1.72) -- (0,0) -- (3.57,1.72) -- cycle    ;
\draw    (150,80) -- (160,90) ;
\draw    (160,90) -- (200,100) ;
\draw    (160,117) -- (160,90) ;
\draw [shift={(160,120)}, rotate = 270] [fill={rgb, 255:red, 0; green, 0; blue, 0 }  ][line width=0.08]  [draw opacity=0] (3.57,-1.72) -- (0,0) -- (3.57,1.72) -- cycle    ;
\draw    (200,100) -- (200,70) ;
\draw    (200,100) -- (210,110) ;
\draw    (210,110) -- (240,110) ;
\draw    (210,127) -- (210,110) ;
\draw [shift={(210,130)}, rotate = 270] [fill={rgb, 255:red, 0; green, 0; blue, 0 }  ][line width=0.08]  [draw opacity=0] (3.57,-1.72) -- (0,0) -- (3.57,1.72) -- cycle    ;
\draw    (240,127) -- (240,110) ;
\draw [shift={(240,130)}, rotate = 270] [fill={rgb, 255:red, 0; green, 0; blue, 0 }  ][line width=0.08]  [draw opacity=0] (3.57,-1.72) -- (0,0) -- (3.57,1.72) -- cycle    ;
\draw    (240,110) -- (247.88,102.12) ;
\draw [shift={(250,100)}, rotate = 135] [fill={rgb, 255:red, 0; green, 0; blue, 0 }  ][line width=0.08]  [draw opacity=0] (3.57,-1.72) -- (0,0) -- (3.57,1.72) -- cycle    ;
\draw    (280,70) -- (320,70) ;
\draw    (320,100) -- (320,70) ;
\draw    (320,70) -- (350,40) ;
\draw    (280,100) -- (320,100) ;
\draw    (320,100) -- (330,110) ;
\draw    (330,110) -- (370,110) ;
\draw    (330,140) -- (330,110) ;
\draw    (370,80) -- (280,80) ;
\draw    (370,80) -- (410,40) ;
\draw    (370,110) -- (370,80) ;
\draw    (370,110) -- (380,120) ;
\draw    (380,140) -- (380,120) ;
\draw    (380,120) -- (420,120) ;
\draw    (420,140) -- (420,120) ;
\draw    (420,120) -- (430,110) ;
\draw    (200,70) -- (207.88,62.12) ;
\draw [shift={(210,60)}, rotate = 135] [fill={rgb, 255:red, 0; green, 0; blue, 0 }  ][line width=0.08]  [draw opacity=0] (3.57,-1.72) -- (0,0) -- (3.57,1.72) -- cycle    ;
\draw  [color={rgb, 255:red, 155; green, 155; blue, 155 }  ,draw opacity=1 ] (280,40) -- (430,40) -- (430,140) -- (280,140) -- cycle ;
\draw    (220,80) -- (258,80) ;
\draw [shift={(260,80)}, rotate = 180] [color={rgb, 255:red, 0; green, 0; blue, 0 }  ][line width=0.75]    (10.93,-3.29) .. controls (6.95,-1.4) and (3.31,-0.3) .. (0,0) .. controls (3.31,0.3) and (6.95,1.4) .. (10.93,3.29)   ;

\end{tikzpicture}

    \caption{A tropical plane curve parametrized by an abstract graph, i.e.\ a tropical stable map.}
    \label{fig-param}
\end{figure}

It turns out that the description of plane tropical curves as subsets of $\mathbb{R}^2$ (with a weighted graph structure) is not sufficient to capture important properties. Rather, it is useful to think of tropical curves as being parametrized by an abstract graph. There can be ambiguity in how to pick a parametrizing graph: in the right picture in Figure \ref{fig-planecurves} for example, one can pick a graph with a $4$-valent vertex and a cycle, or one can pick a tree which is mapped in such a way that two edges cross in the image (i.e. what looks like a $4$-valent vertex in the image is not the image of a vertex, but just the crossing of two edges). This second option to parametrize is depicted in Figure \ref{fig-param}. 
The definition using parametrization by an abstract graph is build in such a way that this ambiguity is resolved.
In the literature, several names are used for these objects, among them \emph{parametrized tropical plane curve} \cite{Mi03} or \emph{tropical stable map} \cite{CJMR17}. The first choice naming is maybe more proverbial, but the second one hints at the close relation to stable maps in Gromov-Witten theory --- a theory which can be applied to the enumerative geometry of plane curves using the idea to represent a plane curve rather as a tuple of an abstract curve and a map to the plane.
In fact, one can view tropical stable maps as degenerations of (log) stable maps as we discuss below, see Remark \ref{rem-tropstablemap}. (We envoke log Gromov-Witten theory to keep track of how curves intersect the coordinate axes of the plane, or, more generally, the toric boundary of a surrounding toric surface \cite{Ran15}.)
We now introduce the parametrizing graphs which are called \emph{abstract tropical curves} and then tropical stable maps. We restrict to the simplest case of rational (i.e. genus 0) curves as we will focus on rational curves in our enumerative study as well.

\begin{definition}[Rational abstract tropical curve]
    A \emph{rational abstract tropical curve} is a metric tree with unbounded edges called ends. The ends have length $\infty$ while all other edges have a length in $\mathbb{R}_{>0}$. We impose a stability condition which requires every vertex to be at least $3$-valent. 
\end{definition}

A \emph{marked rational abstract tropical curve} is a rational abstract tropical curve such that some of its ends are labeled. An isomorphism of marked rational tropical curves is a homeomorphism preserving the lengths of the edges and the markings of ends. 
The \emph{combinatorial type} of a tropical curve is obtained by dropping the information on the metric. 

\begin{definition}[Rational tropical stable map to $\mathbb{R}^2$]\label{def-tropstablemap}
    A \emph{rational tropical stable map to $\mathbb{R}^2$} is a tuple $(\Gamma,f)$ where $\Gamma$ is a marked rational abstract tropical curve and $f:\Gamma\to \mathbb{R}^2$ is a piecewise integer-affine map satisfying:
    \begin{itemize}
        \item On each edge $e$ of $\Gamma$, $f$ is of the form 
        $$t\longmapsto a+t\cdot \direction \mbox{ for some } \direction \in \ZZ^2,$$ 
        where we parametrize $e$ as an interval of size the length $l(e)$ of $e$. The vector $\direction$, called the \emph{direction vector}, arising in this equation is defined up to sign, depending on the starting vertex of the parametrization of the edge. We will sometimes speak of the direction of a flag $\direction(v,e)$, where we consider a flag, i.e.\ a tuple $(v,e)$ of a vertex $v$ and an adjacent edge $e$ as the edge $e$ oriented away from $v$. If $e$ is an end we use the notation $\direction(e)$ for the direction of its unique flag.
        \item The \emph{balancing condition} holds at every vertex, i.e.\ 
        $$\sum_{e \text{ with } v \in \partial e} \direction(v,e)=0.$$
    
    \end{itemize}
\end{definition}
For an edge with direction $\direction=(\mathbf{v}_1,\mathbf{v}_2) \in \ZZ^2$, 
we call $m=\gcd(\mathbf{v}_1,\mathbf{v}_2)$ the \emph{weight} or \emph{expansion factor} and $\frac{1}{m}\cdot \direction$ the \emph{primitive direction} of $e$.

An isomorphism of tropical stable maps is an isomorphism of the underlying tropical curves respecting the maps. 

\begin{example}
    Figure \ref{fig-param} sketches a rational tropical stable map to $\mathbb{R}^2$. The lengths of the bounded edges of the abstract tropical curve are not specified in the picture, but because of Definition \ref{def-tropstablemap} they are uniquely determined by their image.

    There is a unique way to parametrize the left picture of Figure \ref{fig-planecurves} by an abstract tropical curve in such a way that we obtain a tropical stable map to $\mathbb{R}^2$.

\end{example}

\begin{remark}[Drawing convention]
    Even though it is important to work with parametrizations by abstract tropical curves (because otherwise the picture on the right of Figure \ref{fig-planecurves} could be mistaken for a curve of genus $1$), in our pictures, we often just draw the image $f(\Gamma)\subset \mathbb{R}^2$ in the hope that the parametrization is clear just as for the left picture of Figure \ref{fig-planecurves}. 
\end{remark}

\begin{remark}[Dual Newton subdivision]
    Given a tropical stable map $(\Gamma,f)$, the image $f(\Gamma)\subset \mathbb{R}^2$ is  an edge weighted piecewise integer affine linearly embedded graph in $\mathbb{R}^2$ satisfying the balancing condition, and as such dual to a Newton subdivision. The parametrization by an abstract tropical curve cannot be read off from the dual Newton subdivision. For example, the square in the right subdivision in Figure \ref{fig-dual} could be dual to the image of a $4$-valent vertex, or to the crossing of two edges. 
    
    For purposes of enumerative geometry, one can show that we only have to take $3$-valent tropical curves into account (curves with vertices of higher valency only show up for conditions which are not general). We can therefore treat parallelograms in dual Newton subdivisions always as duals to crossings of edges. 
\end{remark}

\begin{remark}[Point conditions]\label{rem-pointconditions}
    In Section \ref{subsec-counttrop}, we introduce tropical enumerative geometry by counting certain rational tropical stable maps satisfying point conditions. Roughly, one can think of a point condition as a point $p\in \mathbb{R}^2$ which is required to lie on the image $f(\Gamma)$ of a tropical stable map $(\Gamma,f)$. Often, it is convenient to keep track of the preimage of $p$ in $\Gamma$. For this purpose, we can use \emph{marked contracted ends} for our abstract tropical curves: these are ends whose direction vector is $\direction(e)=0$. By the balancing condition, assuming the abstract curve is $3$-valent, the two adjacent edges are mapped to opposite directions so that the image locally looks like a tropical curve with a special point, see Figure \ref{fig-point}.
For the purpose of this text, this technical detail can mostly be neglected. 

    \begin{figure}
        \centering

\tikzset{every picture/.style={line width=0.75pt}} 

\begin{tikzpicture}[x=0.75pt,y=0.75pt,yscale=-1,xscale=1]

\draw    (140,180) -- (180,160) ;
\draw  [dash pattern={on 4.5pt off 4.5pt}]  (180,160) -- (180,110) ;
\draw    (180,160) -- (230,170) ;
\draw    (131.79,169.11) -- (170,150) ;
\draw [shift={(130,170)}, rotate = 333.43] [color={rgb, 255:red, 0; green, 0; blue, 0 }  ][line width=0.75]    (6.56,-2.94) .. controls (4.17,-1.38) and (1.99,-0.4) .. (0,0) .. controls (1.99,0.4) and (4.17,1.38) .. (6.56,2.94)   ;
\draw    (238.04,159.61) -- (190,150) ;
\draw [shift={(240,160)}, rotate = 191.31] [color={rgb, 255:red, 0; green, 0; blue, 0 }  ][line width=0.75]    (6.56,-2.94) .. controls (4.17,-1.38) and (1.99,-0.4) .. (0,0) .. controls (1.99,0.4) and (4.17,1.38) .. (6.56,2.94)   ;
\draw    (170,102) -- (170,130) ;
\draw [shift={(170,100)}, rotate = 90] [color={rgb, 255:red, 0; green, 0; blue, 0 }  ][line width=0.75]    (6.56,-2.94) .. controls (4.17,-1.38) and (1.99,-0.4) .. (0,0) .. controls (1.99,0.4) and (4.17,1.38) .. (6.56,2.94)   ;
\draw    (338,150) -- (260,150) ;
\draw [shift={(340,150)}, rotate = 180] [color={rgb, 255:red, 0; green, 0; blue, 0 }  ][line width=0.75]    (10.93,-4.9) .. controls (6.95,-2.3) and (3.31,-0.67) .. (0,0) .. controls (3.31,0.67) and (6.95,2.3) .. (10.93,4.9)   ;
\draw    (380,130) -- (460,150) ;
\draw  [fill={rgb, 255:red, 0; green, 0; blue, 0 }  ,fill opacity=1 ] (418.2,140) .. controls (418.2,139.01) and (419.01,138.2) .. (420,138.2) .. controls (420.99,138.2) and (421.8,139.01) .. (421.8,140) .. controls (421.8,140.99) and (420.99,141.8) .. (420,141.8) .. controls (419.01,141.8) and (418.2,140.99) .. (418.2,140) -- cycle ;
\draw  [color={rgb, 255:red, 155; green, 155; blue, 155 }  ,draw opacity=1 ] (355.2,95.07) -- (492,95.07) -- (492,190.27) -- (355.2,190.27) -- cycle ;
\draw [color={rgb, 255:red, 74; green, 74; blue, 74 }  ,draw opacity=1 ] [dash pattern={on 4.5pt off 4.5pt}]  (365.6,125.87) -- (380,130) ;
\draw [color={rgb, 255:red, 74; green, 74; blue, 74 }  ,draw opacity=1 ] [dash pattern={on 4.5pt off 4.5pt}]  (460,150) -- (474.4,154.13) ;

\draw (144.4,140) node  [font=\footnotesize] [align=left] {\begin{minipage}[lt]{8.67pt}\setlength\topsep{0pt}
$\displaystyle v$
\end{minipage}};
\draw (216,140) node  [font=\footnotesize] [align=left] {\begin{minipage}[lt]{8.67pt}\setlength\topsep{0pt}
$\displaystyle -v$
\end{minipage}};
\draw (156,110) node  [font=\footnotesize] [align=left] {\begin{minipage}[lt]{8.67pt}\setlength\topsep{0pt}
$\displaystyle 0$
\end{minipage}};
\draw (475.2,102.8) node  [font=\footnotesize,color={rgb, 255:red, 155; green, 155; blue, 155 }  ,opacity=1 ] [align=left] {\begin{minipage}[lt]{8.67pt}\setlength\topsep{0pt}
$\displaystyle \mathbb{R}^{2}$
\end{minipage}};
\draw (192,102.8) node  [font=\footnotesize] [align=left] {\begin{minipage}[lt]{8.67pt}\setlength\topsep{0pt}
$\displaystyle x_{1}$
\end{minipage}};
\draw (425.25,127.95) node  [font=\footnotesize] [align=left] {\begin{minipage}[lt]{8.67pt}\setlength\topsep{0pt}
$\displaystyle p_{1}$
\end{minipage}};

\end{tikzpicture}

        \caption{A local picture of a tropical stable map: the marked end $x_1$ meets the point condition $p_1$.}
        \label{fig-point}
    \end{figure}
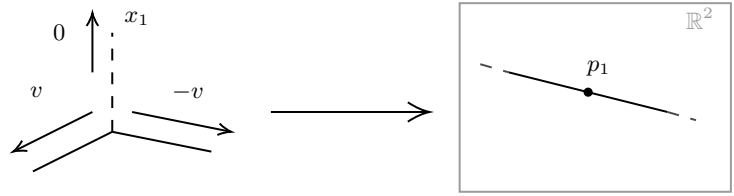
    
    \end{remark}

\begin{definition}[Degree]
The \emph{degree} of a tropical stable map to $\mathbb{R}^2$ is the multiset of direction vectors of its non-contracted ends. 
We restrict to degrees containing the vectors $(-1,0)$, $(0,-1)$ and $(1,1)$ each $d$ times. In this situation, we say that the tropical stable map is of degree $d$.   
Note that this means that all unbounded ends have weight 1.
\end{definition}

\begin{example}\label{ex-degree}
    The two rational tropical stable maps (implicitly) given in Figure \ref{fig-planecurves} are both of degree $3$. This can also be phrased in terms of the dual Newton subdivision: a tropical stable map of degree $d$ yields a subdivision of a triangle with vertices $(0,0)$, $(d,0)$ and $(0,d)$. The fact that such a triangle is the Newton polygon of a general (dehomogenized) polynomial of degree $d$ illuminates the relation.
\end{example}

\begin{remark}[Other degrees]
All results presented in this survey do hold not only for tropical stable maps of degree $d$ (corresponding to curves in $\mathbb{P}^2$ of degree $d$), but for curves in toric del Pezzo surfaces and their tropical counterparts. To ease notation, we restrict to the case of degree $d$ here. 
\end{remark}

\subsection{Tropicalization of plane curves}\label{sec-tropicalization}
Roughly, tropicalization should be viewed as a degeneration process that takes an algebraic object and produces a corresponding tropical object from it. There are many ways in which this rough idea can be made precise. For example, one can consider algebraic curves in $\mathbb{P}^2_{\mathbb{C}}$, remove the points where they meet the coordinate axes, and take a limit of the image under the coordinatewise logarithm of the absolute value map (i.e.\ the limit of a so-called \emph{amoeba}) \cite{Mi03}. To avoid the limit process, one can alternatively study algebraic curves defined over a field with a non-Archimedean valuation:
\begin{definition}[Non-Archimedean valuation]
    Let $K$ be a field. 
    The map $\val:K\rightarrow \mathbb{R}\cup\{\infty\}$ is a \emph{non-Archimedean valuation} if it satisfies:
\begin{align*}
&\val(x)=\infty \Leftrightarrow x=0,\\
&\val(xy)=\val(x)+\val(y),\\
&\val(x+y)\geq \min\{\val(x),\val(y)\}, \mbox{ with equality if } \val(x) \neq \val(y).
\end{align*}
\end{definition}
Our prime example class for fields with a non-Archimedean valuation are fields of Puiseux series with coefficients in a fixed field $k$:

\begin{definition}[Puiseux series]
Let $k$ be a field. The field $\Puiseux{k}$ of \emph{Puiseux series} over $k$ consists of formal series of the form
$$ c(t)=c_1t^{a_1}+c_2 t^{a_2}+\ldots,$$
where the coefficients $c_i$ are in $ k$ and the exponents $a_1<a_2<\ldots \in \QQ$ share a common denominator. (The latter is needed to make the multiplication well-defined.)

The \emph{valuation map} from the field of Puiseux series is:
\begin{align*}\val:&\;\; \Puiseux{k} \longrightarrow \RR\cup \{\infty\}\\
 &\;\;\;c(t)\;\;\longmapsto \;\;a_1 \mbox{ (the least exponent) },\\
&\;\;\;\;\;\;0 \;\;\;\longmapsto  \;\;\infty.
\end{align*}
\end{definition}

We can view tropicalization as taking the valuation componentwise. For conventional reasons, we choose to take minus the valuation rather than the valuation (this corresponds to the so-called max-convention in tropical geometry). The min-convention is also often used, and, for researchers considering certain dualities, it is necessary to work with both conventions \cite{Jos21}. That means that there is no preferred way in the literature. We choose the max-convention since the pictures involving dual Newton subdivisions are a bit clearer. The name max- (resp.\ min-)convention refers to the underlying tropical semifield $(\mathbb{R}\cup\{- \infty\}, \max,+)$ which can be used to set up tropical geometry like algebraic geometry over a ground field and which can be viewed as the closure of the image of $\Puiseux{k}$ of (minus) the valuation map.
Since the tropicalization of zeroes does not play a role for our enumerative purposes in this survey, we restrict tropicalization to the torus points of an algebraic curve:

\begin{definition}[Tropicalization]\label{def-trop}
\emph{Tropicalization} is taking minus the valuation componentwise:
$$\Trop: (\Puiseux{k}^\times)^2\longrightarrow \RR^2: (x,y)\mapsto (-\val(x),-\val(y))$$
where $\Puiseux{k}^\times=\Puiseux{k}\setminus \{0\}$ are the units in $\Puiseux{k}$. 
\end{definition}

\begin{theorem}[Kapranov's theorem, \cite{EKL04}]\label{thm-Kapranov}
Let $C\subset (\Puiseux{k}^\times)^2$ be a curve, where $k$ is algebraically closed of characteristic zero.
Then the closure in the Euclidean topology of the tropicalization 
$\overline {\Trop(C)}$ is a piecewise integer affine linearly embedded graph in $\mathbb{R}^2$ which satisfies the balancing condition with suitable edge weights.
\end{theorem}

\begin{example}[A tropicalized line]
Let $p=x+y+t\in \Puiseux{\CC}[x,y]$. Then $V(p)$ is a line and can be parametrized by
$$V(p)=\big\{(x,-x-t)\;\big|\;x\in \Puiseux{\CC} \big\}.$$
We study the tropicalization of these points.

Case (1):  Let $\val(x)>1$. Then $\val(-x-t)=1$. The points $$\big(-\val(x),-\val(-x-t)\big)$$ of this form yield a ray starting at $(-1,-1)\in \RR^2$ and pointing vertically downwards.

Case (2):  Let $\val(x) <1$. Then $\val(-x-t)=\val(x)$. We obtain a ray starting at $(-1,-1)$ and pointing diagonally left upwards.

Case( 3):  Assume $x=-t+\cdots$. Then $-\val(-x-t)$ can be anything bigger than $1$. We obtain a ray starting at $(-1,-1)$ and pointing leftwards.

Thus, the tropicalization of the line $V(p)$ has a vertex at $(-1,-1)$, and the three rays of direction $(-1,0)$, $(0,-1)$, and $(1,1)$, see Figure \ref{fig-line}.

We can also view it as the image of a tropical stable map of degree $1$.

\begin{figure}
    \centering

\tikzset{every picture/.style={line width=0.75pt}} 

\begin{tikzpicture}[x=0.75pt,y=0.75pt,yscale=-1,xscale=1]

\draw    (160,50) -- (200,50) ;
\draw    (200,90) -- (200,50) ;
\draw    (200,50) -- (230,20) ;

\draw (214,50) node  [font=\footnotesize] [align=left] {\begin{minipage}[lt]{8.67pt}\setlength\topsep{0pt}
$\displaystyle ( -1,-1)$
\end{minipage}};

\end{tikzpicture}

    \caption{The tropicalization of a line is the image of a tropical stable map of degree $1$.}
    \label{fig-line}
\end{figure}
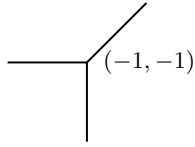

\end{example}

\begin{remark}[Curves over the Puiseux series viewed as families of curves, tropical stable maps as tropicalizations of stable maps]\label{rem-tropstablemap}
     Viewing $t$ as a parameter, a curve $C$ defined over the field of Puiseux series can be considered as a family of curves. Over the reals or complex numbers, one can imagine to insert a small number for t (resp.\ a complex number with a small absolute value). Taking (minus) the valuation then singles out the dominant term of the series and tells us the order of magnitude of the value we can expect when inserting our small number. Carrying this idea further, we can consider a family of curves over a small disc around $0$ and insert a special fiber over $0$. 
     
     Given a curve $C$ defined over the Puiseux series, the theory of semistable degeneration features a reducible nodal curve as special fiber. The generic fiber can be viewed as abstract curve with a morphism to the plane. 
The dual Newton subdivision of the tropicalization of $C$ yields a toric degeneration of the plane to which the special fiber is mapped. More precisely, taking the dual graph of the abstract curve of the special fiber (with a metric taking into account ''how fast'' nodes form, i.e.\ what the valuation of a Puiseux series $p$ in a local equation $xy=p$ of the node is) produces an abstract tropical curve $\Gamma$ \cite{Tyo09}. 
The vertices $v$ of $\Gamma$ correspond to local pieces $C_v$ of the sepcial fiber, defined over a finite field extension of $k$, with a map to the toric surface given by the polygon $\Delta_v$ dual to the image of $v$ in $\mathbb{R}^2$. The bounded edges of $\Gamma$ govern how these local pieces are glued together \cite{Shu06b}. The ends of $\Gamma$ govern the intersection behaviour with the coordinate axes, or, more generally, with the toric boundary of the surrounding toric surface.
In this sense, we can view tropical stable maps as tropicalizations of log stable maps defined over a field with a non-Archimedean valuation such as the Puiseux series over $k$ \cite{Ran15}.
\end{remark}

\section{Complex enumerative geometry and correspondences}\label{sec-complex}
Mikhalkin pioneered tropical geometry as a tool for complex and real enumerative geometry \cite{Mi03}. The prime case to study is the count of rational plane curves of degree $d$ satisfying point conditions, which we introduce in Section \ref{subsec-countcomplex}. In Section \ref{subsec-counttrop}, we introduce the analogous tropical count of curves. Section \ref{sec-tropicalization} already motivated why these counts should be related to each other. In Section \ref{sec-corres} we discuss Mikhalkin's famous correspondence theorem stating the equality of the respective numbers in complex algebraic and tropical geometry.

\subsection{Counting complex plane curves}\label{subsec-countcomplex}
The following enumerative invariants, a priori defined over the complex numbers, are the key players of this survey:

\begin{definition}[The numbers $N_d$]\label{def-Nd}
    Let $N_d$ be the number of rational degree $d$ plane curves subject to $3d-1$ generic point conditions.
\end{definition}
By dimension reasons, $3d-1$ is the right number of points to produce a finite count in this problem rather than a higher-dimensional family. One can show that curves passing through generic point conditions are nodal. For that reason, this count can also be phrased in the context of Gromov-Witten theory and the numbers $N_d$ are also called the Gromov-Witten invariants of the plane. Gromov-Witten theory can also be used to show that these numbers are \emph{invariant}, i.e.\ they do not depend on the chosen point conditions as long as these are generic.

\begin{example}
    The number $N_1$ counts the number of lines through $2$ points, which is $1$. The number $N_2$ counts the number of conics through $5$ generic points, which is also $1$. This can be shown by identifying the space of conics of the form $V(a_0 x^2+a_1 xy+a_2 y^2+a_3 x+a_4 y+a_5)$ with a $\mathbb{P}_\mathbb{C}^5$. Fixing a point condition amounts to taking a codimension one linear subspace in this $\mathbb{P}_{\mathbb{C}}^5$. Fixing five points, we have to take the intersection of five codimension one linear subspaces in $\mathbb{P}_\mathbb{C}^5$, which is one point, corresponding to our one conic through the five points. 

    To determine the number $N_3$, one can take an analogous approach, however, a general cubic curve in $\mathbb{P}_\mathbb{C}^2$ is of genus one. Restricting to rational cubics amounts to restricting to the discriminant hypersurface in the space of all cubics, $\mathbb{P}_\mathbb{C}^9$. The degree of this discriminant hypersurface can be computed to be $12$, and thus, intersecting with $3\cdot 3-1=8$ codimension one linear subspaces for the point conditions again, we obtain $12$ rational plane cubics through $8$ generic points. 

Table \ref{tab-Nd} shows that the numbers $N_d$ grow rather quickly with $d$. Only the first entries of this table can be computed by means of classical algebraic geometry as described above, for the bigger numbers, Kontsevich's famous recursive formula \cite{Kon94} was used which is proved by means of Gromov-Witten theory.

\begin{table}[h!]
  \begin{center}
    \caption{The plane Gromov-Witten invariants $N_d$ for small values of $d$.}
    \label{tab-Nd}
    \begin{tabular}{l|cccccc} $\mbox{degree}$ & $1$&$2$&$3$&$4$&$5$&$6$\\\hline
$N_d$ & $1$ & $1$ & $12$  & $620$ & $87304$ & $ 26312976$\\

    \end{tabular}
  \end{center}
\end{table}
    
\end{example}

In the following, we discuss the analogous tropical counting problem.

\subsection{Counting tropical plane curves}\label{subsec-counttrop}

We essentially repeat Definition \ref{def-Nd} in the tropical world:

\begin{definition}[The tropical count of $N_d$]\label{def-ndtrop}
Fix $3d-1$ generic points in $\mathbb{R}^2$. We define $N_d^{\trop}$ to be the count of rational tropical stable maps $(\Gamma,f)$ to $\mathbb{R}^2$ of degree $d$ that pass through the points (see Remark \ref{rem-pointconditions}), counted with \emph{complex multiplicity}:
$$N_d^\trop= \sum_{(\Gamma,f)} \mult_\CC(\Gamma,f).$$
     \end{definition}

We have to count with multiplicity here, because under tropicalization, several complex curves satisfying point conditions can tropicalize to one and the same tropical curve, see also Section \ref{sec-corres}. If we want our tropical count to reflect the complex count, we have to take this into account. The good news is that this complex multiplicity can be determined purely combinatorially:

\begin{definition}[Complex multiplicity]\label{def-mult}
    Let $(\Gamma,f)$ be a tropical stable degree $d$ map satisfying $3d-1$ generic point conditions and let $v$ be a $3$-valent vertex of $\Gamma$. Assume $x$ and $y$ are the direction vectors of two of the edges adjacent to $v$. Then we define $\mult(v)=|\det(x\ y)|.$
   We define $$\mult_\CC(\Gamma,f)=\prod_v \mult(v),$$ where the product goes over all $3$-valent vertices of $\Gamma$ (which are not adjacent to a contracted end).
\end{definition}
 Because of the balancing condition, it does not matter which two edges adjacent to a vertex $v$ we consider here. In fact, $\mult(v)$ equals the (normalized) area of the triangle dual to the image of $v$ in the dual Newton subdivision.

\begin{example}\label{ex-ninecubics}
Figure \ref{fig-ninecubics} shows $8$ points in $\mathbb{R}^2$ in tropical general position, and the nine rational tropical stable maps of degree $3$ which pass through them. The one in the last row on the right has complex multiplicity $4$, all others have complex multiplicity equal to $1$. We obtain $N^\trop_3= 8\cdot 1 + 1 \cdot 4 =12$ which equals, as expected, the number of complex rational cubics passing through $8$ points. 
 \begin{figure}
\begin{center}
\includegraphics[scale=1.2]{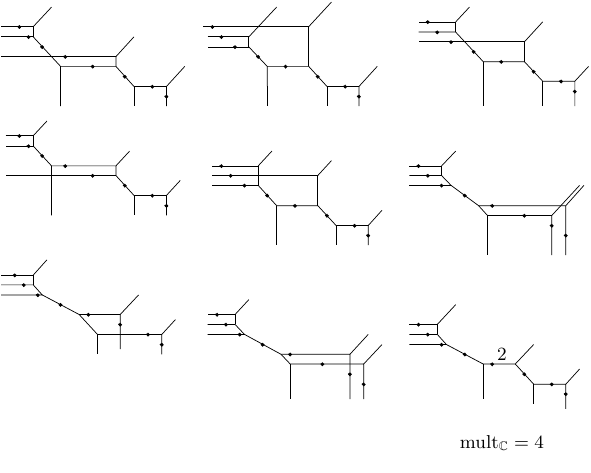}
\end{center}
\caption{The picture shows all 9 rational tropical stable maps of degree $3$ through $8$ given points in general position. The one in the last row on the right has complex multiplicity $4$. All others have complex multiplicity $1$.}\label{fig-ninecubics}
\end{figure}

\end{example}

We still need to show that tropical rational stable maps passing through generic point conditions have only $3$-valent vertices for the definition above to make sense:

\begin{proposition}\label{prop-3valent}
    Let $(\Gamma,f)$ be a tropical rational stable map of degree $d$ passing through $3d-1$ generic points. Then $\Gamma$ is $3$-valent and has no contracted edges (except marked ends). 
\end{proposition}
\begin{proof}[Sketch of proof] 
    $\Gamma$ has $3d$ (uncontracted) ends. We remove the preimages of the point conditions in $\Gamma$, that breaks the tree into $3d$ components. A bounded component would contradict the generic position of the points. Thus every component contains precisely one end. We can orient the edges in each component such that they point away from the point conditions and towards the unique end. The position of the image of a vertex is fixed by two incoming edges. Every vertex must also have one outgoing edge, oriented towards the end. By the genericity of the point conditions, we cannot have more than two incoming edges for a vertex and no contracted bounded edges. 
\end{proof}

\subsection{The correspondence theorem}\label{sec-corres}
Using Definitions \ref{def-Nd} and \ref{def-ndtrop}, we can state the following theorem:

\begin{theorem}[Mikhalkin's Correspondence Theorem \cite{Mi03}]
The complex count of rational degree $d$ plane curves satisfying point conditions coincides with its tropical counterpart:
$$N_d=N_d^\trop.$$    
\end{theorem}

Here, we give a rough overview of the idea of the proof of this theorem.
First, we remark that counting curves over $\mathbb{C}$ or over $\Puiseux{\CC}$ does not make a difference, since the count is the same for any algebraically closed field, and the field of Puiseux series over the complex numbers is algebraically closed.
Now, fix $3d-1$ generic points in $\mathbb{P}^2_{\Puiseux{\CC}}$ and consider the $N_d$ curves that pass through it. Tropicalize them (see Definition \ref{def-trop} and Remark \ref{rem-tropstablemap}). What we obtain is a finite set of tropical rational stable maps that satisfy the tropicalized point conditions, and these are in fact all rational tropical stable maps of degree $d$ that satisfy these point conditions.  Tropicalization then yields a weighted bijection from the set of algebraic curves satisfying the point conditions to the set of tropical rational stable maps satisfying the point conditions. The theorem is proved if we can show that the weights of this weighted bijection are precisely the complex multiplicities of the tropical stable maps. 
To show this, we have to reconstruct all curves that tropicalize to a given tropical stable map $(\Gamma,f)$ \cite{Mi03, Shu04, Shu06b, NS06, AB22, Tyo09, Ran15, Gro16}. For this, we first solve local counting problems for the local pieces $C_v$ of the special fiber of our curve which are mapped to the toric surface dual to the polygon $\Delta_v$ dual to the image of $v$ for a vertex $v$ of $\Gamma$. These local solutions then have to be glued to produce global solutions. The gluing procedure slightly changes in the presence of a point condition.
This reconstruction procedure can also be called \emph{patchworking} of an algebraic curve (even though this phrase was first coined in the real context) \cite{Viro, Shu04, IMS09}.
It is important to note that in this complex setting, we only need to \emph{count solutions} for our local enumerative problems, we do not need to know more about the properties of the local pieces or the gluing, different to the situation in Section \ref{sec-tropW} and \ref{sec-tropQE}.

\section{Welschinger invariants}
With the same tropical techniques, we can also produce real invariant counts of plane curves. In Section \ref{sec-W}, we introduce this real count and discuss its tropical analogue in Section \ref{sec-tropW}.

\subsection{Counting real plane curves}\label{sec-W} 
Real enumerative geometry is in general more difficult than complex enumerative geometry. This is due to the fact that when moving conditions, real solutions may appear or disappear. For example, there can be $8$, $10$ or $12$ real rational cubics passing through $8$ generic real points in the plane \cite{DK00}. However, Welschinger introduced a signed count of real rational curves satisfying real point conditions that is invariant \cite{Wel05} and that is viewed as the real analogue of plane Gromov-Witten invariants. A signed count is of course always a lower bound for the actual number of real curves, which is of interest since one can show using tropical methods that this lower bound is logarithmically equivalent to the Gromov-Witten invariants \cite{IKS04}. 

    A real node of a real curve is called \emph{solitary} if its local equation is $x^2+y^2=0$. 
The local equation of a non-solitary real node is $x^2-y^2=0$, which has two real branches given by $x+y=0$ and $x-y=0$.

\begin{definition}[Welschinger invariant, \cite{Wel05}]
    We define the \emph{Welschinger invariant} to be the signed count of real rational plane curves
    $$W_d = \sum_C (-1)^{m(C)},$$
    where the sum runs over all real curves $C$ of degree $d$ passing through $3d-1$ real generic points. The sign of a curve $C$ is $(-1)^{m(C)}$, where $m(C)$ is the number of real solitary nodes of the curve $C$.
\end{definition}
 This is well-defined because of the following theorem:

\begin{theorem}[\cite{Wel05}]
    The Welschinger invariant $W_d$ is an invariant, i.e.\ it does not depend on the chosen point conditions, as long as they are generic.
\end{theorem}

\subsection{The tropical count of Welschinger invariants}\label{sec-tropW}
What is so nice about the tropical approach to these counting problems is that we can use one and the same set of tropical stable maps no matter whether we discuss real or complex curves:

\begin{definition}\label{def-realmult}
    Let $W_d^\trop=\sum_{(\Gamma,f)}\mult_\RR(\Gamma,f)$, where the sum goes over all tropical rational stable maps of degree $d$ satisfying $3d-1$ generic point conditions, and the real multiplicity is defined to be $$\mult_\RR(\Gamma,f)=\begin{cases}
        0 &\mbox{ if }\mult_\CC (\Gamma,f)\equiv 0 \mod 4,\\
        1 &\mbox{ if }\mult_\CC (\Gamma,f)\equiv 1 \mod 4,\\
        -1 &\mbox{ if }\mult_\CC (\Gamma,f)\equiv 3 \mod 4.\\
    \end{cases}$$
\end{definition}
It is easy to see that the complex multiplicity of a tropical rational stable map of degree $d$ is never congruent to $2\mod 4$: in the dual Newton subdivision, a triangle of an even area must have a boundary edge of even weight by Pick's formula. But since the exterior edges are all of weight $1$, any such triangle is paired with another adjacent even triangle.

\begin{remark}
Note that Definition \ref{def-realmult} appears in other versions in the literature. Firstly, we can express the real multiplicity as a product of vertex multiplicities as in the complex case, where every vertex obtains as real multiplicity $0$, if the complex multiplicity is even, $1$ if it is congruent $1$ modulo $4$, and $-1$ if it is congruent $3$ modulo $4$ (see e.g.\ \cite{Mi03}, Definition 7.19). In other parts of the literature, the real multiplicity of a vertex of odd complex multiplicity is defined to be $(-1)^i$, where $i$ is the number of interior points of the triangle dual to the vertex, see e.g.\ \cite{Shu04}, Proposition 6.1. While these two versions do not coincide vertex by vertex, they coincide when taking the product over all vertices because all ends have weight $1$; this is again due to Pick's formula: The complex multiplicity equals the (normalized) area of the dual triangle which by Pick's formula equals $2i+w_1+w_2+w_3-2$. Here $w_j$ for $j=1,2,3$ denotes the weights of the three adjacent edges. The weight of an edge coincides with the lattice length of the dual edge.
A triangle can have $0$, $1$ or $3$ boundary edges of even lattice length. If there are edges of even lattice length, the real multiplicity is defined to be $0$ in any version. If all triangles have $3$ boundary edges of odd lattice lengths, then $(-1)^{\frac{\mult(v)-1}{2}}=(-1)^{i+\frac{w_1-1}{2}+\frac{w_2-1}{2}+\frac{w_3-1}{2}}$ equals $1$ if $\mult(v)\equiv 1\mod 4$ and $-1$ if  $\mult(v)\equiv 3\mod 4$. As in the total product, non-trivial factors of the form $(-1)^{\frac{w_j-1}{2}}$ can be paired up (since ends have weight $1$, every edge of non-trivial weight is thus bounded and adjacent to $2$ vertices), in the total product over all vertices, the two definitions coincide.
\end{remark}

\begin{example}
    In Figure \ref{fig-ninecubics}, the lower right tropical stable map has real multiplicity $0$ and all others have real multiplicity $1$. We obtain $W_3^\trop=8$.
\end{example}

\begin{theorem}[Mikhalkin's Correspondence Theorem, \cite{Mi03}]
  The Welschinger invariant can be computed tropically, i.e.\
$$W_d=W_d^\trop.$$
\end{theorem}

The idea to prove this correspondence theorem is analogous to the complex version, see Section~\ref{sec-corres}: again, we work over the field of real Puiseux series and establish a weighted bijection between real solutions to our counting problem and the tropical stable maps. To do so, we have to understand the real curves that tropicalize to a given tropical stable map. We do that by constructing local solutions $C_v$ for every vertex $v$ first, and then gluing the local solutions. But this time, it is not sufficient to know \emph{how many} solutions there are --- we must determine how many of our solutions are real, and, furthermore, for each real solution, we must understand how many solitary nodes it has. A careful analysis shows that this can be done vertex by vertex.

\section{Quadratically enriched counts through points defined over the ground field}

\subsection{Background on quadratically enriched counts}\label{sec-QEbasics}

Methods from $\mathbb{A}^1$-homotopy theory allow to define meaningful invariants for enumerative problems posed over a large class of and sometimes even arbitrary base fields $k$. In contrast to the real and complex situation,  instead of producing an integer, one obtains an element of the Grothendieck-Witt ring \(\operatorname{GW}(k)\) of isometry classes of non-degenerate quadratic forms. These quadratic refinements recover the complex invariant via the rank homomorphism \(\operatorname{GW}(k)\to\mathbb{Z}\) and, over \(\mathbb{R}\), recover signed counts via the signature map \(\operatorname{GW}(\mathbb{R})\to\mathbb{Z}\).

We now recall the definition of the Grothendieck-Witt ring $\GW(k)$ in which quadratically enriched counts take their values.  

A \emph{quadratic space} is a finite-dimensional $k$-vector space $V$ equipped with a nondegenerate symmetric bilinear form $q \colon V \times V \to k$. Two quadratic spaces $(V,q)$ and $(V',q')$ are isomorphic if there is an isomorphism of $k$-vector spaces $\phi\colon V \to V'$ such that $q(v,w) = q'(\phi(v), \phi(w))$ for all $v$, $w$ in $V$. 
The set of isomorphism classes of quadratic spaces has the structure of a monoid with respect to direct sum $\oplus$.

\begin{definition}
\label{def:GW}
The  \emph{Grothendieck-Witt ring} $\GW(k)$ is the Grothendieck group of the monoid of quadratic spaces described above. It is a commutative ring with the tensor product as multiplication.
\end{definition}

\begin{definition}
For $a \in k^\times$, we write $\langle a \rangle$ for the class of the $1$-dimensional quadratic space $(k,q)$ with $q(x,y) = axy$.   The {\em hyperbolic plane} is $\mathbb{H} = \langle 1 \rangle + \langle -1 \rangle.$
\end{definition}

\begin{lemma}
\label{lem-GWrelations}
As an additive group, $\GW(k)$ is generated by $\big\{\langle a \rangle : a\in k^\times \big\}$, with relations generated by
\begin{enumerate}
\item $\langle a \rangle= \langle a b^2 \rangle$ for all $a$, $b$ in $k^\times$, and
\item $\langle a\rangle + \langle b \rangle= \langle a+b \rangle+ \langle ab(a+b) \rangle$,
for all $a$, $b$ in $k^\times$ such that $a + b \neq 0$.
\end{enumerate}
Furthermore, the multiplication satisfies the following relation.
\begin{enumerate}
\setcounter{enumi}{2}
\item $\langle a\rangle\langle b\rangle=\langle ab\rangle$ for $a,b\in k^\times$.
\end{enumerate}
\end{lemma}
For a proof, see \cite[Theorem 2.1.11, Remark 2.1.12]{deglise2023notes}.

Assume $k' / k$ is a finite separable field extension. Then any finite dimensional $k'$-vector space $V$ is also finite-dimensional as a $k$-vector space, and we write $V_k$ to denote $V$, viewed as a $k$-vector space. If $(V,q)$ is a quadratic space over $k'$, then $(V_k, \Tr_{k'/k} \circ q)$ is a quadratic space over $k$. One writes 
\[
\Tr_{k'/k} \colon \GW(k') \longrightarrow \GW(k)
\]
for the induced map of Grothendieck-Witt rings.  

\begin{example} \label{ex-GW-CC}
The Grothendieck-Witt ring of the complex numbers is $\GW(\mathbb{C})\simeq \ZZ$ with generator $\langle 1 \rangle$ and the identification is given by evaluating the rank of a quadratic form. The Grothendieck-Witt ring of the reals is $\GW(\mathbb{R})\simeq \ZZ[\mu_2]$, where $\mu_2$ denotes the cyclic group with $2$ elements. It is generated by $\langle 1\rangle$ and $\langle -1\rangle$ where $1=\langle1\rangle$ and $\langle-1\rangle$ is the non-trivial element of $\mu_2$ after the identification with $\mathbb{Z}[\mu_2]$.
\end{example}

One classical result in enumerative geometry is that a smooth cubic surface \(S\subset\mathbb{P}^3\) over an algebraically closed field contains exactly \(27\) lines. 
From the point of view of birational geometry, the cubic surface is rational and it is isomorphic to the blow-up of the projective plane at 6 points in general position (i.e. no 3 points are colinear, and no 6 points lye on a conic) and the 27 lines correspond to the 6 exceptional divisors, the 15 lines that passes through any pair of the 6 points and the 6 conics that pass through any 5 out of the 6 points.

This numerical statement, invariance over algebraically closed fields, does not hold for other fields. A real cubic surface can contain an contain 27, 15, 7 or 3 lines depending on the choice of cubic surface \cite{Schläfli}. In the 40's Segre discovered these real lines can be classified into two different kinds called hyperbolic and elliptic, depending on the dynamics of the intersection points of the line, with the residual conic that appears when intersecting the surfaces with a plane that contains the line under study. Surprisingly, the difference between the number of hyperbolic and elliptic lines equals 3, for every smooth real surface, thus yielding an invariant signed count.

Kass and Wickelgren \cite{KW17} proved that for arbitrary fields that there exists a quadratic weight $\operatorname{ind}_L \in \GW(k(L))$ associated to a line $L$ lying on a smooth cubic surface $S\subset \mathbb{P}^3_k $. They proved that
\[\sum_{L\subset S} \operatorname{Tr}_{k(L)/k} (\operatorname{ind_L})=15\qinv{1}+12\qinv{-1}\in\GW(k),\]
where the sum runs over all scheme-theoretic lines of the surface, or equivalently, over Galois orbits of lines contained in the surface $S$ defined over an algebraic closure of $k$. Such lines are defined over a finite separable field extension $k(L)$ of $k$ and its associated weight $\operatorname{ind}_L\in\GW(k(L))$ is transferred to the base field by the induced trace map $\operatorname{Tr}$.
We remark that this is a non-degenerate quadratic form of rank 27, recovering the natural count for algebraically closed fields. Furthermore, it is of signature 3 for real fields $k\subset\mathbb{R}$. 

\subsection{Quadratically enriched counts of rational plane curves satisfying point conditions defined over the ground field $k$}\label{sec-QE}
From now on we assume that $k$ is a perfect field with $\operatorname{char}k\neq 2,3$.
In contrast to the real and complex situation, one can ask for a weighted count of curves in $\GW(k)$. For example, let $C$ be a rational plane curve defined over $k$ (where $k$ is a perfect field with $\chara(k)\neq 2,3$).  Levine defined the  following quadratic enrichment of the Welschinger sign \cite{LevineWelschinger}: 

\[\operatorname{Wel}_{\mathbb{A}^1}(C)= \Big\langle\prod_{\text{nodes } z}N_{\kappa(z)/\kappa(C)} \big(-\det \operatorname{Hessian}f(z) \big) \Big\rangle \in \GW(\kappa(C))\]
where $\kappa(z)$ is the residue field of the node $z$ and $N_{\kappa(z)/\kappa(C)}$ the field norm.
Note that we take the product over all nodes, not only the ones defined over $\kappa(C)$.
This construction naturally generalizes the Welschinger sign $\operatorname{Wel}(C)$. Indeed, when $\kappa(C)=\mathbb{R}$, a real node $z$ is solitary if and only if $\det\!\operatorname{Hessian}(f)(z)>0$. Since complex nodes contribute a factor of $1$, it follows that
\[
\langle \operatorname{Wel}(C)\rangle=\langle(-1)^{\#\text{real solitary nodes}} \rangle= \operatorname{Wel}_{\A^1}(C)
\]
whenever $\kappa(C)=\mathbb{R}$.

Now define the quadratically enriched count of rational plane curves as 
$$
	N_{\mathbb{A}^1,d}=\sum_{C} \operatorname{Tr}_{\kappa(C)/k}\big(\operatorname{Wel}_{\mathbb{A}^1}(C)\big) \in \GW(k),
$$
where $\kappa(C)$ is the field of definition of $C$ and the sum runs over all $C$ which pass through $3d-1$ generic points defined over $k$. It is shown in \cite{LevineWelschinger, KLSWOrientation} that $N_{\mathbb{A}^1,d}$ is independent of the choice of points.

\begin{remark}[Simultaneous count]
Note that $N_{\mathbb{A}^1,d}$ simultaneously solves our previous counting problems $N_d$ and $W_d$: if we specialize our ground field to $\mathbb{C}$, we take the rank of the quadratic form and obtain $N_d$. If we specialize to $\mathbb{R}$, we take the signature. Our sum defined above also involves complex curves passing through the real points, however, one can show that they always come in pairs being conjugate to each other, and each pair contributes $\mathbb{H}$ to the count of $N_{\mathbb{A}^1,d}$ whose signature is $0$. Thus, when specializing to the reals, it is only the real curves that contribute something non-trivial to the sum, and since they contribute exactly the Welschinger sign, we recover the numbers $W_d$. 

In fact, more is true: when the respective multiplicities for curves are broken down into a product of vertex contributions, then even at that level the rank of the quadratically enriched vertex multiplicity is the complex vertex multiplicity and the signature is the real multiplicity. This implies the statement at the level of tropical stable maps and after summation over all tropical stable maps we obtain the above observation.

\end{remark}

\subsection{The tropical quadratically enriched count and the correspondence theorem}\label{sec-tropQE}
We nearly repeat the first sentence of Section \ref{sec-tropW}: what is so nice about the tropical approach to these counting problems is that we can use the same set of tropical rational stable maps, no matter whether we discuss the real, complex or quadratically enriched count. In fact, the quadratically enriched multiplicity of a tropical stable map is a combination of the complex and real multiplicity:
\begin{definition}
	\label{def-troparith}
Let $(\Gamma,f)$ be a tropical rational stable map of degree $d$ passing through $3d-1$ points in general position.
We define the \emph{quadratically enriched multiplicity of $(\Gamma,f)$}
to be
$$\mult_{\mathbb{A}^1}(\Gamma,f) = \begin{cases} \frac{\mult_{\mathbb{C}}(\Gamma,f)-1}{2}\cdot \mathbb{H}+ \big\langle \mult_{\mathbb{R}}(\Gamma,f) \big\rangle & \mbox{ if $\mult_\CC(\Gamma,f)$ is odd,}\\
\frac{\mult_{\mathbb{C}}(\Gamma,f)}{2}\cdot \mathbb{H} & \mbox{ if $\mult_\CC(\Gamma,f)$ is even.}
\end{cases}$$

We define $$N^\trop_{\mathbb{A}^1,d}=\sum_{(\Gamma,f)} \mult_{\mathbb{A}^1} (\Gamma,f),$$
where the sum goes over all tropical rational stable maps of degree $d$ passing through the $3d-1$ generic points.
\end{definition}

\begin{example}
Using the tropical plane curve count from Example \ref{ex-ninecubics}, Figure \ref{fig-ninecubics}, we obtain $$N^\trop_{\mathbb{A}^1, 3}=2\mathbb{H} + 8\cdot \langle 1 \rangle.$$
\end{example}

\begin{theorem}[Quadratically enriched correspondence theorem, \cite{JPP23}]
\label{thm-corres-new}
Let $k$ be a perfect field of characteristic 0 or large enough, i.e.\ larger than the diameter of $\Delta$.
The correspondence theorem also holds for the quadratically enriched count:
$$N_{\mathbb{A}^1,d}= N^\trop_{\mathbb{A}^1,d}.$$

\end{theorem}

 First let us explain the restriction on the characteristic of the field. To prove the theorem, we break down the multiplicity of a tropical stable map as a product over vertex multiplicities, each involving the weights of the edges. Since $\gw{m}$ is a well-defined element of the Grothendieck-Witt ring only if $m$ is not zero in $k$, we simply restrict the characteristic of $k$ to be larger than the diameter of $\Delta$. 

The proof of this correspondence theorem follows the same outline discussed already in Sections \ref{sec-corres} for the complex version and \ref{sec-tropW} for the real version: we establish a weighted bijection between the set of algebraic curves satisfying the point conditions and the set of tropical rational stable maps satisfying the tropicalized point conditions. We need to show that the weight in this bijection equals the multiplicity of a tropical stable map. To see this, we understand the preimages under tropicalization: little pieces $C_v$ for a vertex $v$ of $\Gamma$, and their gluing. In the complex version, we only needed to know how many such preimages under tropicalization there are. In the real version, we needed to know which preimages are real and how many solitary nodes they have. Now we have to understand the field of definition of every preimage curve, their nodes and their quadratically enriched Welschinger sign. One can see that there is always yet another layer of difficulty for each version of the correspondence theorem!

\section{Algorithms for counting tropical plane curves}\label{sec-alg}

Of course it is great that we can simultaneously solve the counting problems $N_d$, $W_d$ and $N_{\mathbb{A}^1,d}$ with one and the same set of tropical rational stable maps, only adapting their multiplicity to the situation we care for. A central question remains however at this point: how do we actually determine the set of tropical rational stable maps of degree $d$ through $3d-1$ given points?

Algorithms for counting tropical stable maps to $\mathbb{R}^2$ have been studied since the discovery of the correspondence theorem for the numbers $N_d$ \cite{Mi03, Shu06b, BM08, FM09, Blo11, MR08, CJMR17, MR16, CMR25, Gol18}. 
The first such algorithm, Mikhalkin's lattice path algorithm \cite{Mi03}, is based on the idea to establish the possible dual Newton subdivisions of the tropical stable maps passing through our points. Later, a tropical version of the Caporaso-Harris algorithm was studied \cite{GM052, IKS09, JPMPR23}. The original Caporaso-Harris algorithm provides a way to count algebraic plane curves satisfying point conditions by specializing point after point to lie on a line and studying the degenerations of the curves that we obtain from this \cite{CH98}. The tropical version moves one of the point conditions to the very left, which imposes strong restrictions on the combinatorics of the curve, and results in certain splits of tropical stable maps which can be controlled similarly. Iterating the tropical Caporaso-Harris algorithm and distilling the combinatorial essence, one obtains the so-called \emph{floor diagrams} \cite{BM08, FM09, Blo11}.

In this survey, we focus on the lattice path algorithm. It has originally appeared for the complex and real count of tropical stable maps and no quadratically enriched version has appeared yet. Below, we spell out the necessary adaptations, which are obvious to experts, and we propose this quadratically enriched version of the lattice path count as our small original contribution into this survey.

A quadratically enriched version of the Caporaso-Harris algorithm as well as of the quadratically enriched count of floor diagrams has appeared in \cite{JPMPR23}.

\subsection{Quadratically enriched count of lattice paths}
In order to translate the enumeration of tropical stable maps satisfying point conditions to dual Newton subdivisions, we choose a specific (yet general) position for the points. Then, the edges dual to image edges of the tropical stable map which pass through the points form a lattice path. Since we work with tropical stable maps of degree $d$, the dual Newton subdivision is contained in the triangle $\Delta_d$ with vertices $(0,0)$ $(0,d)$ and $(d,0)$ (see Example \ref{ex-degree}).
We order the lattice points of $\Delta_d$ as $(i,j)<(i',j')$ if and only if $i<i'$ or $i=i'$ and $j>j'$.

\begin{definition}[Lattice path]
  A path $\gamma: [0,n] \rightarrow \Delta_d$ for $n \in \NN$ is called a \emph{lattice path}  if
  $\gamma |_{[j-1,j]}$, $j=1,\ldots,n$ is an affine-linear map and $\gamma(j)
  \in\ZZ^2 $ for all $ j=0,\ldots,n $. The path $\gamma$ is called increasing if it respects the order of the lattice points of $\Delta_d$ as discussed above. We require $\gamma(0)=(0,d)$ (the minimal point in our order) and $\gamma(n)=(d,0)$ (the maximal point in our order).
\end{definition}

More generally, one can pick a linear map $\lambda: \RR^2 \to \RR$ whose kernel has an irrational slope, and call a lattice path \emph{$\lambda$-increasing} if $\lambda \circ \gamma$ is strictly increasing. In the definition above, we picked  $\lambda (x,y) = x-\varepsilon y$, where $\varepsilon$ is a small irrational number, then our notion of increasing coincides with $\lambda$-increasing. In this text, we always mean this specific $\lambda$.

\begin{figure}

\begin{center}

\tikzset{every picture/.style={line width=0.75pt}} 

\begin{tikzpicture}[x=0.75pt,y=0.75pt,yscale=-1,xscale=1]

\draw  [fill={rgb, 255:red, 0; green, 0; blue, 0 }  ,fill opacity=1 ] (148.07,50.1) .. controls (148.07,48.9) and (149.04,47.93) .. (150.23,47.93) .. controls (151.43,47.93) and (152.4,48.9) .. (152.4,50.1) .. controls (152.4,51.3) and (151.43,52.27) .. (150.23,52.27) .. controls (149.04,52.27) and (148.07,51.3) .. (148.07,50.1) -- cycle ;
\draw  [fill={rgb, 255:red, 0; green, 0; blue, 0 }  ,fill opacity=1 ] (167.9,70.1) .. controls (167.9,68.9) and (168.87,67.93) .. (170.07,67.93) .. controls (171.26,67.93) and (172.23,68.9) .. (172.23,70.1) .. controls (172.23,71.3) and (171.26,72.27) .. (170.07,72.27) .. controls (168.87,72.27) and (167.9,71.3) .. (167.9,70.1) -- cycle ;
\draw  [fill={rgb, 255:red, 0; green, 0; blue, 0 }  ,fill opacity=1 ] (187.9,90.27) .. controls (187.9,89.07) and (188.87,88.1) .. (190.07,88.1) .. controls (191.26,88.1) and (192.23,89.07) .. (192.23,90.27) .. controls (192.23,91.46) and (191.26,92.43) .. (190.07,92.43) .. controls (188.87,92.43) and (187.9,91.46) .. (187.9,90.27) -- cycle ;
\draw  [fill={rgb, 255:red, 0; green, 0; blue, 0 }  ,fill opacity=1 ] (208.07,110.1) .. controls (208.07,108.9) and (209.04,107.93) .. (210.23,107.93) .. controls (211.43,107.93) and (212.4,108.9) .. (212.4,110.1) .. controls (212.4,111.3) and (211.43,112.27) .. (210.23,112.27) .. controls (209.04,112.27) and (208.07,111.3) .. (208.07,110.1) -- cycle ;
\draw  [fill={rgb, 255:red, 0; green, 0; blue, 0 }  ,fill opacity=1 ] (148.13,70.17) .. controls (148.13,68.97) and (149.1,68) .. (150.3,68) .. controls (151.5,68) and (152.47,68.97) .. (152.47,70.17) .. controls (152.47,71.36) and (151.5,72.33) .. (150.3,72.33) .. controls (149.1,72.33) and (148.13,71.36) .. (148.13,70.17) -- cycle ;
\draw  [fill={rgb, 255:red, 0; green, 0; blue, 0 }  ,fill opacity=1 ] (167.97,90.33) .. controls (167.97,89.14) and (168.94,88.17) .. (170.13,88.17) .. controls (171.33,88.17) and (172.3,89.14) .. (172.3,90.33) .. controls (172.3,91.53) and (171.33,92.5) .. (170.13,92.5) .. controls (168.94,92.5) and (167.97,91.53) .. (167.97,90.33) -- cycle ;
\draw  [fill={rgb, 255:red, 0; green, 0; blue, 0 }  ,fill opacity=1 ] (188.13,110) .. controls (188.13,108.8) and (189.1,107.83) .. (190.3,107.83) .. controls (191.5,107.83) and (192.47,108.8) .. (192.47,110) .. controls (192.47,111.2) and (191.5,112.17) .. (190.3,112.17) .. controls (189.1,112.17) and (188.13,111.2) .. (188.13,110) -- cycle ;
\draw  [fill={rgb, 255:red, 0; green, 0; blue, 0 }  ,fill opacity=1 ] (148.3,90.57) .. controls (148.3,89.37) and (149.27,88.4) .. (150.47,88.4) .. controls (151.66,88.4) and (152.63,89.37) .. (152.63,90.57) .. controls (152.63,91.76) and (151.66,92.73) .. (150.47,92.73) .. controls (149.27,92.73) and (148.3,91.76) .. (148.3,90.57) -- cycle ;
\draw  [fill={rgb, 255:red, 0; green, 0; blue, 0 }  ,fill opacity=1 ] (168.13,109.9) .. controls (168.13,108.7) and (169.1,107.73) .. (170.3,107.73) .. controls (171.5,107.73) and (172.47,108.7) .. (172.47,109.9) .. controls (172.47,111.1) and (171.5,112.07) .. (170.3,112.07) .. controls (169.1,112.07) and (168.13,111.1) .. (168.13,109.9) -- cycle ;
\draw  [fill={rgb, 255:red, 0; green, 0; blue, 0 }  ,fill opacity=1 ] (147.97,110.07) .. controls (147.97,108.87) and (148.94,107.9) .. (150.13,107.9) .. controls (151.33,107.9) and (152.3,108.87) .. (152.3,110.07) .. controls (152.3,111.26) and (151.33,112.23) .. (150.13,112.23) .. controls (148.94,112.23) and (147.97,111.26) .. (147.97,110.07) -- cycle ;
\draw   (30,50) -- (90,110) -- (30,110) -- cycle ;
\draw  [fill={rgb, 255:red, 0; green, 0; blue, 0 }  ,fill opacity=1 ] (247.9,50.22) .. controls (247.9,49.02) and (248.87,48.05) .. (250.07,48.05) .. controls (251.26,48.05) and (252.23,49.02) .. (252.23,50.22) .. controls (252.23,51.41) and (251.26,52.38) .. (250.07,52.38) .. controls (248.87,52.38) and (247.9,51.41) .. (247.9,50.22) -- cycle ;
\draw  [fill={rgb, 255:red, 0; green, 0; blue, 0 }  ,fill opacity=1 ] (267.73,70.22) .. controls (267.73,69.02) and (268.7,68.05) .. (269.9,68.05) .. controls (271.1,68.05) and (272.07,69.02) .. (272.07,70.22) .. controls (272.07,71.41) and (271.1,72.38) .. (269.9,72.38) .. controls (268.7,72.38) and (267.73,71.41) .. (267.73,70.22) -- cycle ;
\draw  [fill={rgb, 255:red, 0; green, 0; blue, 0 }  ,fill opacity=1 ] (287.73,90.38) .. controls (287.73,89.19) and (288.7,88.22) .. (289.9,88.22) .. controls (291.1,88.22) and (292.07,89.19) .. (292.07,90.38) .. controls (292.07,91.58) and (291.1,92.55) .. (289.9,92.55) .. controls (288.7,92.55) and (287.73,91.58) .. (287.73,90.38) -- cycle ;
\draw  [fill={rgb, 255:red, 0; green, 0; blue, 0 }  ,fill opacity=1 ] (307.9,110.22) .. controls (307.9,109.02) and (308.87,108.05) .. (310.07,108.05) .. controls (311.26,108.05) and (312.23,109.02) .. (312.23,110.22) .. controls (312.23,111.41) and (311.26,112.38) .. (310.07,112.38) .. controls (308.87,112.38) and (307.9,111.41) .. (307.9,110.22) -- cycle ;
\draw  [fill={rgb, 255:red, 0; green, 0; blue, 0 }  ,fill opacity=1 ] (247.97,70.28) .. controls (247.97,69.09) and (248.94,68.12) .. (250.13,68.12) .. controls (251.33,68.12) and (252.3,69.09) .. (252.3,70.28) .. controls (252.3,71.48) and (251.33,72.45) .. (250.13,72.45) .. controls (248.94,72.45) and (247.97,71.48) .. (247.97,70.28) -- cycle ;
\draw  [fill={rgb, 255:red, 0; green, 0; blue, 0 }  ,fill opacity=1 ] (267.8,90.45) .. controls (267.8,89.25) and (268.77,88.28) .. (269.97,88.28) .. controls (271.16,88.28) and (272.13,89.25) .. (272.13,90.45) .. controls (272.13,91.65) and (271.16,92.62) .. (269.97,92.62) .. controls (268.77,92.62) and (267.8,91.65) .. (267.8,90.45) -- cycle ;
\draw  [fill={rgb, 255:red, 0; green, 0; blue, 0 }  ,fill opacity=1 ] (287.97,110.12) .. controls (287.97,108.92) and (288.94,107.95) .. (290.13,107.95) .. controls (291.33,107.95) and (292.3,108.92) .. (292.3,110.12) .. controls (292.3,111.31) and (291.33,112.28) .. (290.13,112.28) .. controls (288.94,112.28) and (287.97,111.31) .. (287.97,110.12) -- cycle ;
\draw  [fill={rgb, 255:red, 0; green, 0; blue, 0 }  ,fill opacity=1 ] (248.13,90.68) .. controls (248.13,89.49) and (249.1,88.52) .. (250.3,88.52) .. controls (251.5,88.52) and (252.47,89.49) .. (252.47,90.68) .. controls (252.47,91.88) and (251.5,92.85) .. (250.3,92.85) .. controls (249.1,92.85) and (248.13,91.88) .. (248.13,90.68) -- cycle ;
\draw  [fill={rgb, 255:red, 0; green, 0; blue, 0 }  ,fill opacity=1 ] (267.97,110.02) .. controls (267.97,108.82) and (268.94,107.85) .. (270.13,107.85) .. controls (271.33,107.85) and (272.3,108.82) .. (272.3,110.02) .. controls (272.3,111.21) and (271.33,112.18) .. (270.13,112.18) .. controls (268.94,112.18) and (267.97,111.21) .. (267.97,110.02) -- cycle ;
\draw  [fill={rgb, 255:red, 0; green, 0; blue, 0 }  ,fill opacity=1 ] (247.8,110.18) .. controls (247.8,108.99) and (248.77,108.02) .. (249.97,108.02) .. controls (251.16,108.02) and (252.13,108.99) .. (252.13,110.18) .. controls (252.13,111.38) and (251.16,112.35) .. (249.97,112.35) .. controls (248.77,112.35) and (247.8,111.38) .. (247.8,110.18) -- cycle ;
\draw  [fill={rgb, 255:red, 0; green, 0; blue, 0 }  ,fill opacity=1 ] (348.33,50.05) .. controls (348.33,48.85) and (349.3,47.88) .. (350.5,47.88) .. controls (351.7,47.88) and (352.67,48.85) .. (352.67,50.05) .. controls (352.67,51.25) and (351.7,52.22) .. (350.5,52.22) .. controls (349.3,52.22) and (348.33,51.25) .. (348.33,50.05) -- cycle ;
\draw  [fill={rgb, 255:red, 0; green, 0; blue, 0 }  ,fill opacity=1 ] (368.17,70.05) .. controls (368.17,68.85) and (369.14,67.88) .. (370.33,67.88) .. controls (371.53,67.88) and (372.5,68.85) .. (372.5,70.05) .. controls (372.5,71.25) and (371.53,72.22) .. (370.33,72.22) .. controls (369.14,72.22) and (368.17,71.25) .. (368.17,70.05) -- cycle ;
\draw  [fill={rgb, 255:red, 0; green, 0; blue, 0 }  ,fill opacity=1 ] (388.17,90.22) .. controls (388.17,89.02) and (389.14,88.05) .. (390.33,88.05) .. controls (391.53,88.05) and (392.5,89.02) .. (392.5,90.22) .. controls (392.5,91.41) and (391.53,92.38) .. (390.33,92.38) .. controls (389.14,92.38) and (388.17,91.41) .. (388.17,90.22) -- cycle ;
\draw  [fill={rgb, 255:red, 0; green, 0; blue, 0 }  ,fill opacity=1 ] (408.33,110.05) .. controls (408.33,108.85) and (409.3,107.88) .. (410.5,107.88) .. controls (411.7,107.88) and (412.67,108.85) .. (412.67,110.05) .. controls (412.67,111.25) and (411.7,112.22) .. (410.5,112.22) .. controls (409.3,112.22) and (408.33,111.25) .. (408.33,110.05) -- cycle ;
\draw  [fill={rgb, 255:red, 0; green, 0; blue, 0 }  ,fill opacity=1 ] (348.4,70.12) .. controls (348.4,68.92) and (349.37,67.95) .. (350.57,67.95) .. controls (351.76,67.95) and (352.73,68.92) .. (352.73,70.12) .. controls (352.73,71.31) and (351.76,72.28) .. (350.57,72.28) .. controls (349.37,72.28) and (348.4,71.31) .. (348.4,70.12) -- cycle ;
\draw  [fill={rgb, 255:red, 0; green, 0; blue, 0 }  ,fill opacity=1 ] (368.23,90.28) .. controls (368.23,89.09) and (369.2,88.12) .. (370.4,88.12) .. controls (371.6,88.12) and (372.57,89.09) .. (372.57,90.28) .. controls (372.57,91.48) and (371.6,92.45) .. (370.4,92.45) .. controls (369.2,92.45) and (368.23,91.48) .. (368.23,90.28) -- cycle ;
\draw  [fill={rgb, 255:red, 0; green, 0; blue, 0 }  ,fill opacity=1 ] (388.4,109.95) .. controls (388.4,108.75) and (389.37,107.78) .. (390.57,107.78) .. controls (391.76,107.78) and (392.73,108.75) .. (392.73,109.95) .. controls (392.73,111.15) and (391.76,112.12) .. (390.57,112.12) .. controls (389.37,112.12) and (388.4,111.15) .. (388.4,109.95) -- cycle ;
\draw  [fill={rgb, 255:red, 0; green, 0; blue, 0 }  ,fill opacity=1 ] (348.57,90.52) .. controls (348.57,89.32) and (349.54,88.35) .. (350.73,88.35) .. controls (351.93,88.35) and (352.9,89.32) .. (352.9,90.52) .. controls (352.9,91.71) and (351.93,92.68) .. (350.73,92.68) .. controls (349.54,92.68) and (348.57,91.71) .. (348.57,90.52) -- cycle ;
\draw  [fill={rgb, 255:red, 0; green, 0; blue, 0 }  ,fill opacity=1 ] (368.4,109.85) .. controls (368.4,108.65) and (369.37,107.68) .. (370.57,107.68) .. controls (371.76,107.68) and (372.73,108.65) .. (372.73,109.85) .. controls (372.73,111.05) and (371.76,112.02) .. (370.57,112.02) .. controls (369.37,112.02) and (368.4,111.05) .. (368.4,109.85) -- cycle ;
\draw  [fill={rgb, 255:red, 0; green, 0; blue, 0 }  ,fill opacity=1 ] (348.23,110.02) .. controls (348.23,108.82) and (349.2,107.85) .. (350.4,107.85) .. controls (351.6,107.85) and (352.57,108.82) .. (352.57,110.02) .. controls (352.57,111.21) and (351.6,112.18) .. (350.4,112.18) .. controls (349.2,112.18) and (348.23,111.21) .. (348.23,110.02) -- cycle ;
\draw    (250.07,50.22) -- (249.97,110.18) ;
\draw    (250,110) -- (310,110) ;
\draw    (350,50) -- (410,110) ;

\draw (135.17,40.05) node [anchor=north west][inner sep=0.75pt]  [font=\scriptsize] [align=left] {$\displaystyle 1$};
\draw (135.17,61.22) node [anchor=north west][inner sep=0.75pt]  [font=\scriptsize] [align=left] {$\displaystyle 2$};
\draw (136,81.55) node [anchor=north west][inner sep=0.75pt]  [font=\scriptsize] [align=left] {$\displaystyle 3$};
\draw (135,101.38) node [anchor=north west][inner sep=0.75pt]  [font=\scriptsize] [align=left] {$\displaystyle 4$};
\draw (158.83,61.38) node [anchor=north west][inner sep=0.75pt]  [font=\scriptsize] [align=left] {$\displaystyle 5$};
\draw (159,81.22) node [anchor=north west][inner sep=0.75pt]  [font=\scriptsize] [align=left] {$\displaystyle 6$};
\draw (158.5,101.22) node [anchor=north west][inner sep=0.75pt]  [font=\scriptsize] [align=left] {$\displaystyle 7$};
\draw (178.67,81.72) node [anchor=north west][inner sep=0.75pt]  [font=\scriptsize] [align=left] {$\displaystyle 8$};
\draw (178.33,101.22) node [anchor=north west][inner sep=0.75pt]  [font=\scriptsize] [align=left] {$\displaystyle 9$};
\draw (194.17,101.72) node [anchor=north west][inner sep=0.75pt]  [font=\scriptsize] [align=left] {$\displaystyle 10$};
\draw (231,117) node [anchor=north west][inner sep=0.75pt]  [font=\scriptsize] [align=left] {$\displaystyle \delta _{-} :[ 0,6]\rightarrow \Delta _{d}$};
\draw (361,52) node [anchor=north west][inner sep=0.75pt]  [font=\scriptsize] [align=left] {$\displaystyle \delta _{+} :[ 0,3]\rightarrow \Delta _{d}$};

\end{tikzpicture}

\end{center}

 \caption{The triangle $\Delta_3$, the order of its lattice points, and the two paths $\delta_-$ and $\delta_+$ on its boundary.}
    \label{fig:endpaths}
\end{figure}
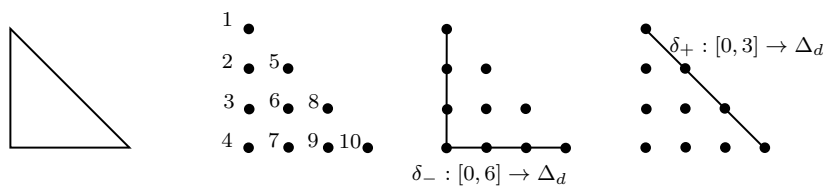

 The minimal and maximal points of $\Delta_d$, that is $(0,d)$ and $(d,0)$, divide the boundary $\partial \Delta_d$ into two
increasing lattice paths $\delta_{+}:[0,d]\rightarrow \partial
\Delta_d$ (going clockwise around $\partial \Delta_d$) and
$\delta_-:[0,2d]\rightarrow \partial \Delta_d$ (going counterclockwise around
$\partial \Delta_d$). 
Figure \ref{fig:endpaths} shows an example for the triangle $\Delta_3$. 

To make use of dual Newton subdivisions in our count, we need to express the quadratically enriched multiplicity of a tropical stable map as a product of vertex multiplicities. 

\begin{definition}[Quadratically enriched vertex multiplicity]
Let $v$ be a vertex of a tropical stable map and assume it has $3$ adjacent edges of weights $w_1$, $w_2$, $w_3$.
We define the \emph{quadratically enriched multiplicity of $v$} to be
$$\mult_{\mathbb{A}^1}(v) = 
\begin{cases} 
	\frac{\mult(v)-1}{2}\cdot \mathbb{H}+ \big\langle (-1)^i \cdot w_1 \cdot w_2\cdot w_3 \big\rangle & \mbox{ if $\mult(v)$ is odd,}\\
	\frac{\mult(v)}{2}\cdot \mathbb{H} & \mbox{ if $\mult(v)$ is even,}
\end{cases}$$
where $\mult(v)$ denotes the complex multiplicity of the vertex, see Definition \ref{def-mult} and $i$ denotes the number of interior lattice points in the triangle $\Delta_v$ dual to the image of $v$.
Since the definition only depends on the triangle $\Delta_v$ dual to the image of $v$, we also use the notation $\mult_{\mathbb{A}^1}(\Delta_v)$.
\end{definition}

\begin{lemma}[\cite{JPMPR23}, Lemma 3.3]\label{lem-a1mult}
    Let $(\Gamma,f)$ be a tropical rational stable map of degree $d$ passing through $3d-1$ generic points. Then $\mult_{\mathbb{A}^1}(\Gamma,f)=\prod_v\mult_{\mathbb{A}^1}(v),$
    where the product goes over all $3$-valent vertices $v$ of $\Gamma$ (which are not adjacent to a contracted end).
\end{lemma}

Since the quadratically enriched multiplicity is expressed in terms of the complex and real multiplicity and the equality above is clear when we specialize to the complex numbers, to prove \cref{lem-a1mult} one only needs to understand the real multiplicity in terms of vertex multiplicities, see \cite{Mi03, IKS09}.

\begin{definition} \label{def-mu}
  Let $\gamma:[0,n]\rightarrow \Delta_d$ be an increasing path from the minimal to the maximal point. The \emph{(positive and negative) quadratically enriched multiplicities} 
  $\mu_+(\gamma)$ and $\mu_-(\gamma)$ are defined recursively as follows:
  \begin{enumerate}
  \item \label{def-mu-a}
    $\mu_{\pm}(\delta_{\pm}):=\langle 1 \rangle$.
  \item \label{def-mu-b}
    If $\gamma \neq \delta_{\pm}$, let $k_{\pm} \in [0,n]$ be the smallest
    number such that $\gamma$ makes a left turn (respectively a right turn for
    $\mu_-$) at $\gamma(k_{\pm})$. (If no such $k_{\pm}$ exists we set
    $\mu_{\pm}(\gamma):= 0$). We define two other $\lambda$-increasing lattice
    paths $\gamma'$ and $\gamma''$ as follows (compare Figure \ref{fig:twopaths}):
    \begin {itemize}
    \item $\gamma_{\pm}':[0,n-1]\rightarrow \Delta$ is the path that cuts the
      corner of $\gamma(k_{\pm})$, i.e.\ $\gamma'_{\pm}(j):=\gamma(j)$ for
      $j<k_{\pm}$ and $\gamma'_{\pm}(j) := \gamma(j+1)$ for $j \geq k_{\pm}$.
    \item $\gamma''_{\pm}:[0,n]\rightarrow \Delta$ is the path that completes
      the corner of $\gamma(k_{\pm})$ to a parallelogram, i.e.\ $\gamma''_{\pm}
      (j):=\gamma(j)$ for all $j \neq k_{\pm}$ and $\gamma''_{\pm}(k_{\pm})
      :=\gamma(k_{\pm}-1)+\gamma(k_{\pm}+1)-\gamma(k_{\pm})$.
      \end{itemize}

      \begin{figure}         
            \begin{center}
            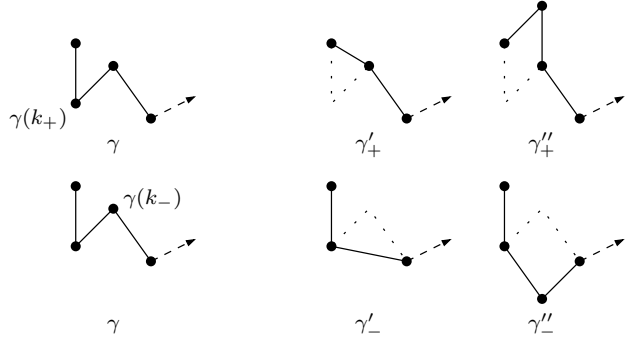
            \end{center}
 
          \caption{The two $\lambda$-increasing paths $\gamma'$ and $\gamma''$ formed recursively from $\gamma$.}
          \label{fig:twopaths}
      \end{figure}
      
   \noindent
    Let $\Delta_v$ be the triangle with vertices $\gamma(k_{\pm}-1),\gamma(k_{\pm}),
    \gamma(k_{\pm}+1)$. Then we set
    \begin{equation*}
      \mu_{\pm}(\gamma):=\mult_{\mathbb{A}^1}(\Delta_v) \cdot \mu_{\pm}(\gamma'_{\pm})
        +\mu_{\pm}(\gamma''_{\pm}).
    \end{equation*}
    As both paths $\gamma'_{\pm}$ and $\gamma''_{\pm}$ include a smaller area
    with $\delta_{\pm}$, we can assume that their multiplicity is known. If
    $\gamma''_{\pm}$ does not map to $\Delta$, then $\mu_{\pm}(\gamma''_{\pm})$ is
    defined to be 0.
  \end{enumerate}
  Finally, the \emph{quadratically enriched multiplicity} $\mu(\gamma)$ is defined to be the product
  $\mu(\gamma) := \mu_+(\gamma) \mu_-(\gamma)$.
\end{definition}
Note that if the recursion terminates with a path which is neither $\delta_+$ nor $\delta_-$, then this counts zero, see Figure \ref{fig:endpaths}. Paths without a left (resp.\ right) turn, but not equal to $\delta_{\pm}$, or ``faster'' paths (taking fewer but bigger steps)
$\delta':[0,d']\rightarrow \Delta$ such that
$\delta_+([0,d])=\delta'([0,d'])$ but $d'<d$ have multiplicity zero.
\begin{remark}\label{rem-suitnewton}
One can interpret the recursion for computing the multiplicity of a path in terms of Newton subdivisions of $\Delta$. Let $\gamma$ be an increasing path. First, consider the positive multiplicity $\mu_+$. The two paths $\gamma$ and $\gamma'$ enclose the triangle $\Delta_v$. The two paths $\gamma$ and $\gamma''$ enclose a parallelogram. Take two copies of the Newton polygon $\Delta$ and draw the triangle in one of them and the parallelogram in the other. Iterating this process, we can fill in different choices for dual subdivisions consisting of triangles and parallelograms of the part of $\Delta$ above $\gamma$. Analogously, one obtains subdivisions of the part of $\Delta$ below $\gamma$ when performing the recursion for $\mu_-$. Any such subdivision above $\gamma$ can be combined with any subdivision below $\gamma$. The set of subdivisions which arise like this is called the set of \emph{possible Newton subdivisions for $\gamma$}. The multiplicity $\mu$ of a path $\gamma$ is nothing else but the number of possible Newton subdivisions for $\gamma$ counted with a multiplicity that equals the quadratically enriched multiplicity of a tropical stable map to which this subdivision is dual.
\end{remark}
\begin{example}
 Figure \ref{fig:possiblesubd} shows an increasing path $\gamma$ and the two possible
  Newton subdivisions for $\gamma$. Both possible subdivisions are dual to tropical stable maps of quadratically enriched multiplicity
  $\langle 1 \rangle$, they are depicted in the left column of Figure \ref{fig-ninecubics}, on the top and middle.
  \begin{figure}

\begin{center}

\tikzset{every picture/.style={line width=0.75pt}} 

\begin{tikzpicture}[x=0.75pt,y=0.75pt,yscale=-1,xscale=1]

\draw  [fill={rgb, 255:red, 0; green, 0; blue, 0 }  ,fill opacity=1 ] (247.9,50.22) .. controls (247.9,49.02) and (248.87,48.05) .. (250.07,48.05) .. controls (251.26,48.05) and (252.23,49.02) .. (252.23,50.22) .. controls (252.23,51.41) and (251.26,52.38) .. (250.07,52.38) .. controls (248.87,52.38) and (247.9,51.41) .. (247.9,50.22) -- cycle ;
\draw  [fill={rgb, 255:red, 0; green, 0; blue, 0 }  ,fill opacity=1 ] (267.73,70.22) .. controls (267.73,69.02) and (268.7,68.05) .. (269.9,68.05) .. controls (271.1,68.05) and (272.07,69.02) .. (272.07,70.22) .. controls (272.07,71.41) and (271.1,72.38) .. (269.9,72.38) .. controls (268.7,72.38) and (267.73,71.41) .. (267.73,70.22) -- cycle ;
\draw  [fill={rgb, 255:red, 0; green, 0; blue, 0 }  ,fill opacity=1 ] (287.73,90.38) .. controls (287.73,89.19) and (288.7,88.22) .. (289.9,88.22) .. controls (291.1,88.22) and (292.07,89.19) .. (292.07,90.38) .. controls (292.07,91.58) and (291.1,92.55) .. (289.9,92.55) .. controls (288.7,92.55) and (287.73,91.58) .. (287.73,90.38) -- cycle ;
\draw  [fill={rgb, 255:red, 0; green, 0; blue, 0 }  ,fill opacity=1 ] (307.9,110.22) .. controls (307.9,109.02) and (308.87,108.05) .. (310.07,108.05) .. controls (311.26,108.05) and (312.23,109.02) .. (312.23,110.22) .. controls (312.23,111.41) and (311.26,112.38) .. (310.07,112.38) .. controls (308.87,112.38) and (307.9,111.41) .. (307.9,110.22) -- cycle ;
\draw  [fill={rgb, 255:red, 0; green, 0; blue, 0 }  ,fill opacity=1 ] (247.97,70.28) .. controls (247.97,69.09) and (248.94,68.12) .. (250.13,68.12) .. controls (251.33,68.12) and (252.3,69.09) .. (252.3,70.28) .. controls (252.3,71.48) and (251.33,72.45) .. (250.13,72.45) .. controls (248.94,72.45) and (247.97,71.48) .. (247.97,70.28) -- cycle ;
\draw  [fill={rgb, 255:red, 0; green, 0; blue, 0 }  ,fill opacity=1 ] (267.8,90.45) .. controls (267.8,89.25) and (268.77,88.28) .. (269.97,88.28) .. controls (271.16,88.28) and (272.13,89.25) .. (272.13,90.45) .. controls (272.13,91.65) and (271.16,92.62) .. (269.97,92.62) .. controls (268.77,92.62) and (267.8,91.65) .. (267.8,90.45) -- cycle ;
\draw  [fill={rgb, 255:red, 0; green, 0; blue, 0 }  ,fill opacity=1 ] (287.97,110.12) .. controls (287.97,108.92) and (288.94,107.95) .. (290.13,107.95) .. controls (291.33,107.95) and (292.3,108.92) .. (292.3,110.12) .. controls (292.3,111.31) and (291.33,112.28) .. (290.13,112.28) .. controls (288.94,112.28) and (287.97,111.31) .. (287.97,110.12) -- cycle ;
\draw  [fill={rgb, 255:red, 0; green, 0; blue, 0 }  ,fill opacity=1 ] (248.13,90.68) .. controls (248.13,89.49) and (249.1,88.52) .. (250.3,88.52) .. controls (251.5,88.52) and (252.47,89.49) .. (252.47,90.68) .. controls (252.47,91.88) and (251.5,92.85) .. (250.3,92.85) .. controls (249.1,92.85) and (248.13,91.88) .. (248.13,90.68) -- cycle ;
\draw  [fill={rgb, 255:red, 0; green, 0; blue, 0 }  ,fill opacity=1 ] (267.97,110.02) .. controls (267.97,108.82) and (268.94,107.85) .. (270.13,107.85) .. controls (271.33,107.85) and (272.3,108.82) .. (272.3,110.02) .. controls (272.3,111.21) and (271.33,112.18) .. (270.13,112.18) .. controls (268.94,112.18) and (267.97,111.21) .. (267.97,110.02) -- cycle ;
\draw  [fill={rgb, 255:red, 0; green, 0; blue, 0 }  ,fill opacity=1 ] (247.83,110) .. controls (247.83,108.8) and (248.8,107.83) .. (250,107.83) .. controls (251.2,107.83) and (252.17,108.8) .. (252.17,110) .. controls (252.17,111.2) and (251.2,112.17) .. (250,112.17) .. controls (248.8,112.17) and (247.83,111.2) .. (247.83,110) -- cycle ;
\draw    (250,50) -- (250.3,90.68) ;
\draw    (269.9,70.22) -- (250.3,90.68) ;
\draw    (270,110) -- (269.9,70.22) ;
\draw    (289.9,90.22) -- (270.3,110.68) ;
\draw    (290,110) -- (289.9,90.22) ;
\draw    (290,110) -- (309.9,110.22) ;
\draw  [fill={rgb, 255:red, 0; green, 0; blue, 0 }  ,fill opacity=1 ] (488.56,50.06) .. controls (488.56,48.86) and (489.53,47.89) .. (490.72,47.89) .. controls (491.92,47.89) and (492.89,48.86) .. (492.89,50.06) .. controls (492.89,51.25) and (491.92,52.22) .. (490.72,52.22) .. controls (489.53,52.22) and (488.56,51.25) .. (488.56,50.06) -- cycle ;
\draw  [fill={rgb, 255:red, 0; green, 0; blue, 0 }  ,fill opacity=1 ] (508.39,70.06) .. controls (508.39,68.86) and (509.36,67.89) .. (510.56,67.89) .. controls (511.75,67.89) and (512.72,68.86) .. (512.72,70.06) .. controls (512.72,71.25) and (511.75,72.22) .. (510.56,72.22) .. controls (509.36,72.22) and (508.39,71.25) .. (508.39,70.06) -- cycle ;
\draw  [fill={rgb, 255:red, 0; green, 0; blue, 0 }  ,fill opacity=1 ] (528.39,90.22) .. controls (528.39,89.03) and (529.36,88.06) .. (530.56,88.06) .. controls (531.75,88.06) and (532.72,89.03) .. (532.72,90.22) .. controls (532.72,91.42) and (531.75,92.39) .. (530.56,92.39) .. controls (529.36,92.39) and (528.39,91.42) .. (528.39,90.22) -- cycle ;
\draw  [fill={rgb, 255:red, 0; green, 0; blue, 0 }  ,fill opacity=1 ] (548.56,110.06) .. controls (548.56,108.86) and (549.53,107.89) .. (550.72,107.89) .. controls (551.92,107.89) and (552.89,108.86) .. (552.89,110.06) .. controls (552.89,111.25) and (551.92,112.22) .. (550.72,112.22) .. controls (549.53,112.22) and (548.56,111.25) .. (548.56,110.06) -- cycle ;
\draw  [fill={rgb, 255:red, 0; green, 0; blue, 0 }  ,fill opacity=1 ] (488.62,70.12) .. controls (488.62,68.93) and (489.59,67.96) .. (490.79,67.96) .. controls (491.99,67.96) and (492.96,68.93) .. (492.96,70.12) .. controls (492.96,71.32) and (491.99,72.29) .. (490.79,72.29) .. controls (489.59,72.29) and (488.62,71.32) .. (488.62,70.12) -- cycle ;
\draw  [fill={rgb, 255:red, 0; green, 0; blue, 0 }  ,fill opacity=1 ] (508.46,90.29) .. controls (508.46,89.09) and (509.43,88.12) .. (510.62,88.12) .. controls (511.82,88.12) and (512.79,89.09) .. (512.79,90.29) .. controls (512.79,91.49) and (511.82,92.46) .. (510.62,92.46) .. controls (509.43,92.46) and (508.46,91.49) .. (508.46,90.29) -- cycle ;
\draw  [fill={rgb, 255:red, 0; green, 0; blue, 0 }  ,fill opacity=1 ] (528.62,109.96) .. controls (528.62,108.76) and (529.59,107.79) .. (530.79,107.79) .. controls (531.99,107.79) and (532.96,108.76) .. (532.96,109.96) .. controls (532.96,111.15) and (531.99,112.12) .. (530.79,112.12) .. controls (529.59,112.12) and (528.62,111.15) .. (528.62,109.96) -- cycle ;
\draw  [fill={rgb, 255:red, 0; green, 0; blue, 0 }  ,fill opacity=1 ] (488.79,90.52) .. controls (488.79,89.33) and (489.76,88.36) .. (490.96,88.36) .. controls (492.15,88.36) and (493.12,89.33) .. (493.12,90.52) .. controls (493.12,91.72) and (492.15,92.69) .. (490.96,92.69) .. controls (489.76,92.69) and (488.79,91.72) .. (488.79,90.52) -- cycle ;
\draw  [fill={rgb, 255:red, 0; green, 0; blue, 0 }  ,fill opacity=1 ] (508.62,109.86) .. controls (508.62,108.66) and (509.59,107.69) .. (510.79,107.69) .. controls (511.99,107.69) and (512.96,108.66) .. (512.96,109.86) .. controls (512.96,111.05) and (511.99,112.02) .. (510.79,112.02) .. controls (509.59,112.02) and (508.62,111.05) .. (508.62,109.86) -- cycle ;
\draw  [fill={rgb, 255:red, 0; green, 0; blue, 0 }  ,fill opacity=1 ] (488.49,109.84) .. controls (488.49,108.64) and (489.46,107.67) .. (490.66,107.67) .. controls (491.85,107.67) and (492.82,108.64) .. (492.82,109.84) .. controls (492.82,111.04) and (491.85,112.01) .. (490.66,112.01) .. controls (489.46,112.01) and (488.49,111.04) .. (488.49,109.84) -- cycle ;
\draw    (490.66,49.84) -- (490.96,90.52) ;
\draw    (510.56,70.06) -- (490.96,90.52) ;
\draw    (510.66,109.84) -- (510.56,70.06) ;
\draw    (530.56,90.06) -- (510.96,110.52) ;
\draw    (530.66,109.84) -- (530.56,90.06) ;
\draw    (530.66,109.84) -- (550.56,110.06) ;
\draw  [fill={rgb, 255:red, 0; green, 0; blue, 0 }  ,fill opacity=1 ] (407.84,50.39) .. controls (407.84,49.19) and (408.81,48.22) .. (410.01,48.22) .. controls (411.21,48.22) and (412.18,49.19) .. (412.18,50.39) .. controls (412.18,51.59) and (411.21,52.56) .. (410.01,52.56) .. controls (408.81,52.56) and (407.84,51.59) .. (407.84,50.39) -- cycle ;
\draw  [fill={rgb, 255:red, 0; green, 0; blue, 0 }  ,fill opacity=1 ] (427.68,70.39) .. controls (427.68,69.19) and (428.65,68.22) .. (429.84,68.22) .. controls (431.04,68.22) and (432.01,69.19) .. (432.01,70.39) .. controls (432.01,71.59) and (431.04,72.56) .. (429.84,72.56) .. controls (428.65,72.56) and (427.68,71.59) .. (427.68,70.39) -- cycle ;
\draw  [fill={rgb, 255:red, 0; green, 0; blue, 0 }  ,fill opacity=1 ] (447.68,90.56) .. controls (447.68,89.36) and (448.65,88.39) .. (449.84,88.39) .. controls (451.04,88.39) and (452.01,89.36) .. (452.01,90.56) .. controls (452.01,91.75) and (451.04,92.72) .. (449.84,92.72) .. controls (448.65,92.72) and (447.68,91.75) .. (447.68,90.56) -- cycle ;
\draw  [fill={rgb, 255:red, 0; green, 0; blue, 0 }  ,fill opacity=1 ] (467.84,110.39) .. controls (467.84,109.19) and (468.81,108.22) .. (470.01,108.22) .. controls (471.21,108.22) and (472.18,109.19) .. (472.18,110.39) .. controls (472.18,111.59) and (471.21,112.56) .. (470.01,112.56) .. controls (468.81,112.56) and (467.84,111.59) .. (467.84,110.39) -- cycle ;
\draw  [fill={rgb, 255:red, 0; green, 0; blue, 0 }  ,fill opacity=1 ] (407.91,70.46) .. controls (407.91,69.26) and (408.88,68.29) .. (410.08,68.29) .. controls (411.27,68.29) and (412.24,69.26) .. (412.24,70.46) .. controls (412.24,71.65) and (411.27,72.62) .. (410.08,72.62) .. controls (408.88,72.62) and (407.91,71.65) .. (407.91,70.46) -- cycle ;
\draw  [fill={rgb, 255:red, 0; green, 0; blue, 0 }  ,fill opacity=1 ] (427.74,90.62) .. controls (427.74,89.43) and (428.71,88.46) .. (429.91,88.46) .. controls (431.11,88.46) and (432.08,89.43) .. (432.08,90.62) .. controls (432.08,91.82) and (431.11,92.79) .. (429.91,92.79) .. controls (428.71,92.79) and (427.74,91.82) .. (427.74,90.62) -- cycle ;
\draw  [fill={rgb, 255:red, 0; green, 0; blue, 0 }  ,fill opacity=1 ] (447.91,110.29) .. controls (447.91,109.09) and (448.88,108.12) .. (450.08,108.12) .. controls (451.27,108.12) and (452.24,109.09) .. (452.24,110.29) .. controls (452.24,111.49) and (451.27,112.46) .. (450.08,112.46) .. controls (448.88,112.46) and (447.91,111.49) .. (447.91,110.29) -- cycle ;
\draw  [fill={rgb, 255:red, 0; green, 0; blue, 0 }  ,fill opacity=1 ] (408.08,90.86) .. controls (408.08,89.66) and (409.05,88.69) .. (410.24,88.69) .. controls (411.44,88.69) and (412.41,89.66) .. (412.41,90.86) .. controls (412.41,92.05) and (411.44,93.02) .. (410.24,93.02) .. controls (409.05,93.02) and (408.08,92.05) .. (408.08,90.86) -- cycle ;
\draw  [fill={rgb, 255:red, 0; green, 0; blue, 0 }  ,fill opacity=1 ] (427.91,110.19) .. controls (427.91,108.99) and (428.88,108.02) .. (430.08,108.02) .. controls (431.27,108.02) and (432.24,108.99) .. (432.24,110.19) .. controls (432.24,111.39) and (431.27,112.36) .. (430.08,112.36) .. controls (428.88,112.36) and (427.91,111.39) .. (427.91,110.19) -- cycle ;
\draw  [fill={rgb, 255:red, 0; green, 0; blue, 0 }  ,fill opacity=1 ] (407.78,110.17) .. controls (407.78,108.98) and (408.75,108.01) .. (409.94,108.01) .. controls (411.14,108.01) and (412.11,108.98) .. (412.11,110.17) .. controls (412.11,111.37) and (411.14,112.34) .. (409.94,112.34) .. controls (408.75,112.34) and (407.78,111.37) .. (407.78,110.17) -- cycle ;
\draw    (409.94,50.17) -- (409.94,110.17) ;
\draw    (429.84,70.39) -- (410.24,90.86) ;
\draw    (429.94,110.17) -- (429.84,70.39) ;
\draw    (449.84,90.39) -- (430.24,110.86) ;
\draw    (449.94,110.17) -- (449.84,90.39) ;
\draw    (449.94,110.17) -- (469.84,110.39) ;
\draw    (409.94,70.17) -- (429.84,70.39) ;
\draw    (429.94,90.17) -- (449.84,90.39) ;
\draw    (490.79,70.12) -- (510.69,70.34) ;
\draw    (510.79,90.12) -- (530.69,90.34) ;
\draw    (409.94,50.17) -- (470.01,110.39) ;
\draw    (490.72,50.06) -- (550.79,110.27) ;
\draw    (429.94,90.17) -- (409.94,110.17) ;
\draw    (429.94,110.17) -- (449.84,110.39) ;
\draw    (410.18,109.97) -- (430.08,110.19) ;
\draw    (490.72,90.07) -- (510.62,90.29) ;
\draw    (490.66,109.84) -- (510.56,110.06) ;
\draw    (510.89,109.74) -- (530.79,109.96) ;
\draw    (490.66,109.84) -- (490.96,90.52) ;

\end{tikzpicture}

  \end{center}
 \caption{A lattice path $\gamma$ and its two possible Newton subdivisions. They are dual to the two tropical stable maps depicted in the left column on the top and middle in Figure \ref{fig-ninecubics}.}
      \label{fig:possiblesubd}
  \end{figure}

\end{example}
\begin{definition}\label{def-npath}
   We denote by $ N_{\mathbb{A}^1,d}^{\paths}  $ the
  \emph{ count of increasing lattice paths} $\gamma:[0,3d-1]
  \rightarrow \Delta_d$ counted with
  their quadratically enriched  multiplicities, as in Definition \ref {def-mu}.

\end{definition}

\begin{theorem} \label{thm-correstroppath}
The quadratically enriched count of lattice paths in $\Delta_d$ equals the quadratically enriched count of tropical rational stable maps of degree $d$ through $3d-1$ generic points:
  $$N_{\mathbb{A}^1,d}^{\paths}=N_{\mathbb{A}^1,d}^{\trop},$$ where $N_{\mathbb{A}^1,d}^{\paths}$ is defined in Definition~\ref{def-npath} and $N^{\trop}_{\mathbb{A}^1,d}$ is defined in Definition~\ref{def-troparith}.
\end{theorem}
The complex version of this statement was originally proved in \cite{Mi03}, Theorem 2.

The proof establishes a bijection between possible Newton subdivisions for each increasing path (see Remark \ref{rem-suitnewton}) and tropical stable maps passing through a special choice of $3d-1$ generic points, such that the bijection respects the multiplicity. The bijection is the same, no matter whether we speak about a complex, real, or quadratically enriched count. Our contribution is only to adapt the multiplicity of a lattice path to the quadratically enriched setting.
In our outline of the proof, we follow the outline of the complex version of the proof in \cite{CMR23}, Chapter 10. Some of the pictures are taken from there.

The following definition concerns the special choice of point conditions needed for the proof of Theorem \ref{thm-correstroppath}.
\begin{definition}\label{def-plambda}
Choose a line $H$ orthogonal to the
  kernel of $\lambda:\RR^2\rightarrow \RR: (x,y)\mapsto x-\varepsilon y$ and $3d-1$ generic points $p_1,\ldots,p_{3d-1}$ on $H$ such
  that the distance between $p_i$ and $p_{i+1}$ is much bigger than the
  distance between $p_{i-1}$ and $p_i$ for all $i$. A set of points  $\mathcal{P}_{\lambda} = \{p_1,\ldots,p_{3d-1}\}$ with this property is said to be in \emph{Mikhalkin position}.
\end{definition}

\begin{lemma}\label{lem-H}
Let $(\Gamma,f)$ be a tropical stable map through $\mathcal{P}_{\lambda}$. Then $f(\Gamma)$ intersects the line $H$ (on which the points $\mathcal{P}_{\lambda}=(p_1,\ldots,p_{3d-1})$ lie) only in the points $p_1,\ldots,p_{3d-1}$.
\end{lemma}
This can be shown using the fact that every component of $\Gamma$ minus the preimages of the point conditions has a unique end, which was already used in Proposition \ref{prop-3valent}. For more details, see \cite{CMR23}, Lemma 10.1.9.

\begin{lemma}\label{lem-curvepath}
Let $(\Gamma,f)$ be a tropical rational stable map of degree $d$ through $3d-1$ points in Mikhalkin position $\mathcal{P}_{\lambda}$. 
Let $\Xi $ denote the edges in the dual Newton subdivision of $f(\Gamma)$ which are dual to the edges passing through $\mathcal{P}_{\lambda}$.
Then $\Xi$ is the image of an increasing path $\gamma:[0,3d-1]\rightarrow \Delta_d$.
\end{lemma}
\begin{proof}
First note that by Lemma \ref{lem-H}, the set $\mathcal{P}_{\lambda}$ coincides with $f(\Gamma)\cap H$. Thus one can equivalently show that the edges in the subdivision which are dual to edges of $f(\Gamma)$ which intersect $H$ form a lattice path.

Consider a vertex $V$ of the Newton subdivision and the edges adjacent to it. Dual to these edges is a chain of edges of $f(\Gamma)$ which encloses a convex polyhedron. Any line meets this chain of edges at most twice. 
We distinguish several cases depending on the position of $V$.
\begin{enumerate}
\item \label{case1} Assume first $V$ is in the interior of $\Delta$. 
\begin{figure}[tb]
\begin{center}
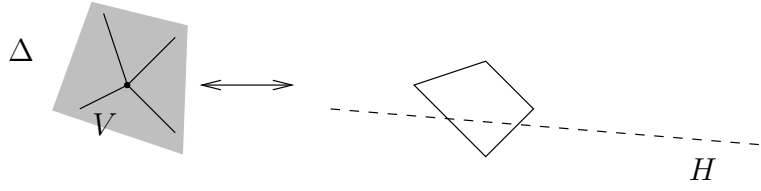

\end{center}
 \caption{Case (\ref{case1}) in the proof of Lemma \ref{lem-curvepath}.}
    \label{fig:case1}
\end{figure}
Then the convex polyhedron is bounded and $H$ meets either none or two of the dual edges, see Figure \ref{fig:case1}. It cannot meet a vertex, as $\mathcal{P}_{\lambda}$ is in general position. Hence, either none or two edges of $\Xi$ must be adjacent to $V$.

\item \label{case2} Assume next that $V=(0,d)$.
Recall that $(0,d)$ is the vertex of $\Delta_d$ where $\lambda$ attains its minimum. If we draw a line parallel to the kernel of $\lambda$ through $(0,d)$ then this line meets $\Delta_d$ only in $(0,d)$. 
Even more, the edges adjacent to $(0,d)$ in $\Delta_d$ lie on one side of the line parallel to $\ker \lambda$.
Assume $H$ intersects two edges of $f(\Gamma)$ which are adjacent to $(0,d)$. 
Change the coordinate system for a moment such that $H$ is of slope $0$. Then the slope of one edge must be negative and the slope of the other edge must be positive, see Figure \ref{fig:case2a}.

\begin{figure}
   
\begin{center}

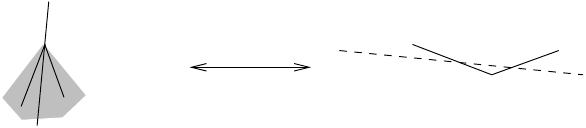

\end{center}

    \caption{Case (\ref{case2}a) of the proof of Lemma \ref{lem-curvepath}.}
    \label{fig:case2a}
\end{figure}

But then the duals of these edges in  $\Delta_d$ would not be on one side of a line parallel to the kernel of $\lambda$.
So it is not possible that $H$ intersects more than one of the dual edges to the edges adjacent to $(0,d)$.

Assume $H$ intersects none of the dual edges to the interior edges adjacent to $(0,d)$. Then either all those edges lie above $H$, or below $H$. Without restriction, assume they lie above $H$. Then also the ends dual to the edges in the boundary of $\Delta_d$ adjacent to $(0,d)$ lie above $H$. Also, as both of these edges lie on one side of the line parallel to $\ker\lambda$, one of the dual ends has to intersect $H$, see Figure \ref{fig:case2b}.
\begin{figure}
   
\begin{center}
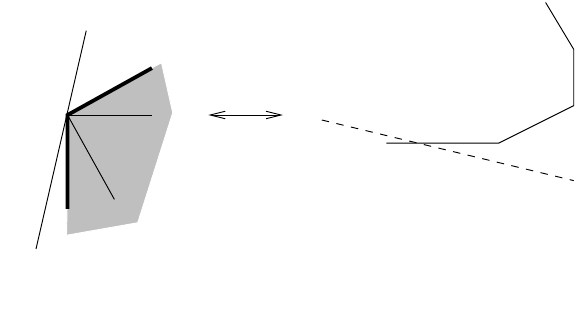

\end{center}

    \caption{Case (\ref{case2}b) of the proof of Lemma \ref{lem-curvepath}.}
    \label{fig:case2b}
\end{figure}

Hence altogether $H$ intersects precisely one of the duals of the edges adjacent to $(0,d)$. So there is one edge of $\Xi$ adjacent to $(0,d)$.
\item Assume that $V=(d,0)$, then we get the same result as for $(0,d)$: there is one edge of $\Xi$ adjacent to $(d,0)$.
\item Assume finally that $V$ is in the boundary of $\Delta_d$, but neither equal to $(0,d)$ nor $(d,0)$.

Then the edges adjacent to $V$ do not lie on one side of a line parallel to $\ker \lambda$ through $V$. Therefore we can see as before that $H$ intersects two of these edges, if it intersects any at all.
So there are either two or no edge of $\Xi$ adjacent to $V$.
\end{enumerate}
Altogether this shows that at each vertex $V$ --- except $(0,d)$ and $(d,0)$ --- there are either two  edges of $\Xi$ adjacent or none, while at $(0,d)$ and $(d,0)$, there is one edge adjacent.  

Finally, we show that $\Xi$ is increasing, and thus a path from $(0, d)$ to $(d, 0)$. Assume the vertices $a_1$, $a_2$ and $a_3$ are three consecutive vertices of $\Xi$, such that the step from $a_1$ to $a_2$ is increasing, while the step from $a_2$ to $a_3$ is not. But this means, that $a_1$ and $a_3$ lie on the same side of a line parallel to $\ker \lambda$ through $a_2$. By the above, $H$ cannot intersect both dual edges, which contradicts the assumption that the two edges were part of $\Xi$.
\end{proof}

\begin {proof}[Proof of Theorem \ref{thm-correstroppath}:]
Take a point configuration $\mathcal{P}_{\lambda}$ in Mikhalkin position (see Definition \ref{def-plambda}). Consider the tropical rational stable maps $(\Gamma,f)$ of degree $d$ that pass through these points.
If we take the edges of $f(\Gamma)$ that pass through $\mathcal{P}_{\lambda}$ and consider their dual edges in the Newton subdivision then these
  dual edges form an increasing path by Lemma \ref{lem-curvepath}.

Let $\gamma$ be a path. The following shows that there are exactly $\mu(\gamma)$ tropical stable maps (counted with quadratically enriched multiplicity) through $\mathcal{P}_{\lambda}$, such that the edges passing through the points are dual to $\Xi=\gamma$. 

 Interpret $\image(\gamma)$ as a set of edges in the subdivision and try to draw a dual tropical stable map. For the edges passing through the points of $\mathcal{P}_{\lambda}$, the direction is prescribed by the path $\gamma$. For each, a point through which it should pass is prescribed by $\mathcal{P}_{\lambda}$. 
Take the first edge of $\Xi=\image(\gamma)$ --- that is, the one starting at $(0,d)$ --- and draw a (part of a) line orthogonal to this marked edge through $p_1$. Going on, draw a line dual to the next edge of $\Xi=\image(\gamma)$ through $p_2$ and so on. Treat these line segments as the germs of  edges of the image of a tropical stable map dual to a possible subdivision for $\Xi$. Recursively, the edge germs grow into $f(\Gamma)$ dual to a possible subdivision for $\Xi$. To prove the result, one has to count the possibilities for this. There will be one tropical stable map for each possible Newton subdivision for $\gamma$ (see Remark \ref{rem-suitnewton}). 

The multiplicity $\mu_+$
  counts the possibilities to complete the edges dual to $\gamma$ to $f(\Gamma)$ (weighted with their quadratically enriched multiplicity) in the
  half-plane above $H$, whereas $\mu_-$ counts below $H$. 
The following makes this argument precise for $\mu_+$ (for $\mu_-$ it is analogous).

Let the first left turn of the path $\gamma$ be enclosed by the edges $E$ and $E'$ whose duals $e$ and $e'$ pass through $p_i$ and $p_{i+1}$. The edges through the points $p_1,\ldots,p_{i-1}$ do not intersect above $H$, as this was the first left turn. The edges $e$ and $e'$  intersect above $H$, but below all other possible intersections of dual edges of $\Xi=\image(\gamma)$. This is true due to the chosen configuration of points: the distance of $p_{j+1}$ and $p_j$ is much bigger than the distance of $p_j$ and $p_{j-1}$. That is, one can draw a parallel line $H'$ to $H$ such that $H$ and $H'$ enclose a strip in which only the intersection point of $e$ and $e'$ lies. 
Passing from $\gamma$
  to $\gamma'$ and $\gamma''$ corresponds to moving the line $H$ up to $H'$. The path $\gamma'$ leaves a triangle $\Delta_v$ out, and $\gamma''$ completes the corner to a parallelogram. These two possibilities are  dual to the two possibilities how $f(\Gamma)$ can look like at $e\cap e'$: it can either have a $3$-valent vertex $v$ --- in which case it is dual to the triangle $\Delta_v$ --- or $e$ and $e'$ can just intersect --- in which case it is dual to the parallelogram which is enclosed by $\gamma$ and $\gamma''$.
So the change from $\gamma$ to $\gamma'$ and $\gamma''$ describes the possibilities how the image of a tropical stable map through $\mathcal{P}_{\lambda}$ can look like in the strip enclosed by $H$ and $H'$.

Figure \ref{fig:strip} shows a path $\gamma$ and the two paths $\gamma'$ and $\gamma''$, together with the triangle, respectively, parallelogram which they enclose with $\gamma$. Below, the dual curves in the strip enclosed by $H$ and $H'$ are shown. 

\begin{figure}
  
\begin{center}
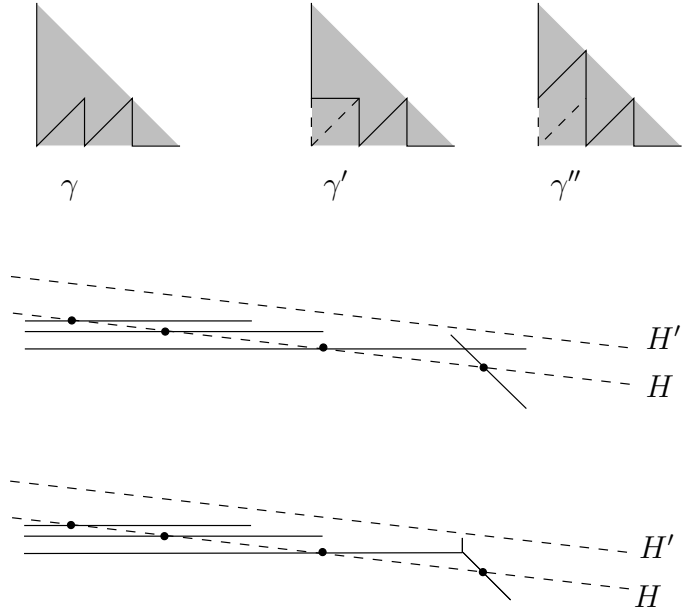

\end{center}

    \caption{The possibilities for the image of a tropical stable map in the strip formed by the line $H$ and a parallel line correspond to the two paths $\gamma'$ and $\gamma''$ in the recursion to compute the quadratically enriched multiplicity of a path, see Definition \ref{def-mu}.}
    \label{fig:strip}
\end{figure}

Recursively, we can see that there is in fact exactly one tropical stable map  through $\mathcal{P}_{\lambda}$ whose image $f(\Gamma)$ is dual to each possible Newton subdivision for $\gamma$. This curve $(\Gamma,f)$ is obviously of degree $d$, as we end up with the two paths $\delta_+$ and $\delta_-$ which are dual to the ends prescribed by $\Delta_d$ in our recursion. 

Hence we constructed a (weighted) bijection between the set of possible Newton subdivisions for an increasing path and the set of tropical rational stable maps of degree $d$ through $\mathcal{P}_{\lambda}$. The quadratically enriched multiplicity of a tropical rational stable map can be expressed in terms of the dual Newton subdivision as $\prod \mult_{\mathbb{A}^1}(\Delta_v)$ by Lemma \ref{lem-a1mult}. By the recursive definition of the quadratically enriched multiplicity of a lattice path, we can understand it as a sum over all possible Newton subdivisions for this path, where each such Newton subdivision contributes the multiplicity $\prod \mult_{\mathbb{A}^1}(\Delta_v)$. Thus the bijection respects the quadratically enriched multiplicities of tropical rational stable maps resp.\ lattice paths and the statement follows.
\end {proof}
\begin{example}
Figure \ref{fig:expaths} shows an increasing path $\gamma:[0,8]\rightarrow \Delta_3$. 
Next to it, its two possible Newton subdivisions are shown, and below the dual curves through the point configuration $\mathcal{P}_{\lambda}$ in Mikhalkin position. The distance between the points $p_{i+1}$ and $p_i$ is not very much bigger than the distance between $p_i$ and $p_{i-1}$ in the picture, just because this is hard to depict. The picture should still be sufficient to give an idea on how the point configuration in Mikhalkin position and the dual tropical stable maps for the possible subdivisions for our paths that pass through the points look like. A more accurate picture, for which the distances of the points grows more, can be produced from this by a topological deformation, stretching horizontally.
\begin{figure}
   
\begin{center}
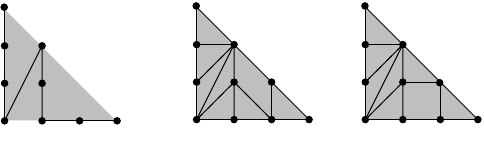
\end{center}
\begin{center}
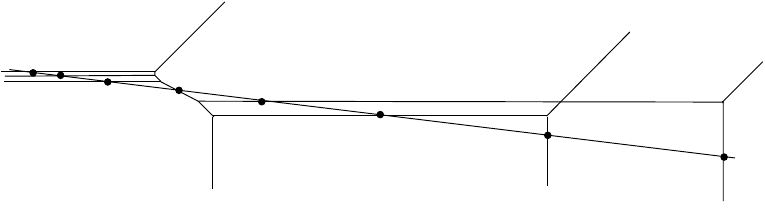
\end{center}
\begin{center}
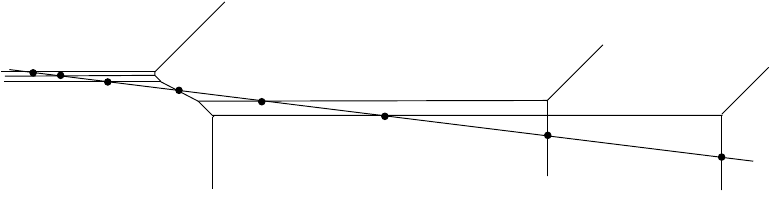
\end{center}

    \caption{An increasing path, its two possible dual Newton subdivisions, and their corresponding tropical plane curves passing through points in Mikhalkin position.}
    \label{fig:expaths}
\end{figure}

\end{example}

\section{Quadratically enriched counts of curves satisfying point conditions defined over quadratic extensions of the ground field}\label{sec-QEQ}

\subsection{The algebraic geometry version of the count and its background}
In the real world, it is very natural to count real curves that satisfy real point conditions, as the signed count for $W_d$ described in Section \ref{sec-W}. If a real plane curve $V(f)$, with homogeneous $f\in \RR[x,y,z]$, passes through a point $p$ defined over the complex numbers, i.e.\ $f(p)=0$, then also $f(\overline{p})=0$ for the complex conjugate point $\overline{p}$. From that point of view, it also makes sense to study counts of rational degree $d$ real curves satisfying $r$ real points and $s$ pairs of complex conjugate point conditions, where $r+2s=3d-1$. If we count these curves with the same sign as before, $(-1)^m$, where $m$ is the number of solitary nodes, we obtain again a Welschinger invariant, denoted by $W_{d,s}$, so that $W_{d,0}=W_d$.
In \cite{Shu06b}, Shustin developed a tropical correspondence and a tropical way to compute the numbers $W_{d,s}$, by means of a generalized lattice path algorithm.

When introducing quadratically enriched counts, we first focused on counts of curves satisfying point conditions which are defined over the ground field.
Similarly to the extension of Welschinger invariants, also this quadratically enriched count has been extended recently, by Kass, Levine, Solomon and Wickelgren \cite{KassLevineSolomonWickelgren}, see also \cite{LevineWelschinger}:
Let $\sigma=(L_1,\ldots,L_m)$ be the list of finite field extensions of the ground field $k$ over which the point conditions are defined.
Then  $N_{\AA^1,d}(\sigma)$ denotes the count of rational plane curves of degree $d$ satisfying these point conditions, with quadratically enriched multiplicity as in Section \ref{sec-QE}.
Setting $k = \RR$ and $\sigma=(\RR^r,\CC^s)$, taking the signature of $N_{\AA^1,d}(\sigma)$, one recovers the Welschinger invariant $W_{d,s}$, and taking the rank one recovers the number $N_d$. 

$N_{\AA^1,d}(\sigma)$ is an enumerative invariant, i.e.\ does not depend on the chosen point conditions, as long as they are generic and as prescribed by $\sigma$. Developing effective means for computing these invariants still poses challenges. Given how nicely the tropical approach fits to the quadratically enriched counts of curves satisfying $k$-points (see Section \ref{sec-tropQE}), it is natural to try the tropical approach, as we will describe in Section \ref{sec-tropQEQ}. Other possibilities for computing these invariants have been studied in \cite{BW25, BRW25}.

\subsection{The tropical version of the quadratically enriched counts of curves satisfying point conditions defined over quadratic extensions of the ground field}\label{sec-tropQEQ}
To make the tropical approach feasible, we first restrict to the case of quadratic field extensions for our point conditions, see also Remark \ref{rem-otherfieldextensions}.

If we consider a point $p$ defined over $\Puiseux{\CC}$ and its conjugate $\overline{p}$, then they obviously tropicalize to the same point in $\mathbb{R}^2$. The same holds true for conjugate points in quadratic field extensions. For this reason, we have to consider \emph{double} point conditions for our tropical stable maps $(\Gamma,f)$ now: point conditions which have to be met by two points (i.e.\ contracted ends) of $\Gamma$.

While the exact description of the features of a tropical stable map satisfying such double point conditions is combinatorially quite involved, we feel that the basic idea to just merge two point conditions and make the tropical stable maps follow along is sufficiently intuitive. Hence, we refrain from giving a precise definition of the tropical stable maps to be counted and only present an example. For more details, we refer to \cite{JPMPR25}, Section 3 for the discussion of the combinatorics underlying tropical stable maps satisfying double point conditions, and to Section 4 for the definition of the quadratically enriched multiplicity of these tropical stable maps. Theorem 1.2 in \cite{JPMPR25} states the equality of the thus defined tropical numbers with the invariants  $N_{\AA^1,d}(\sigma)$ and Theorem~1.4 provides a tropical algorithm build on floor diagrams to compute them.

\begin{figure}
\begin{center}
\tikzset{every picture/.style={line width=1pt}} 

\begin{tikzpicture}[x=0.75pt,y=0.75pt,yscale=-2.5,xscale=2.5]

\draw  [fill={rgb, 255:red, 0; green, 0; blue, 0 }  ,fill opacity=1 ] (122.49,154.14) .. controls (122.49,153.36) and (123.1,152.72) .. (123.85,152.72) .. controls (124.6,152.72) and (125.22,153.36) .. (125.22,154.14) .. controls (125.22,154.92) and (124.6,155.56) .. (123.85,155.56) .. controls (123.1,155.56) and (122.49,154.92) .. (122.49,154.14) -- cycle ;
\draw  [fill={rgb, 255:red, 0; green, 0; blue, 0 }  ,fill opacity=1 ] (179.94,168.23) .. controls (179.94,167.45) and (180.55,166.81) .. (181.3,166.81) .. controls (182.06,166.81) and (182.67,167.45) .. (182.67,168.23) .. controls (182.67,169.01) and (182.06,169.65) .. (181.3,169.65) .. controls (180.55,169.65) and (179.94,169.01) .. (179.94,168.23) -- cycle ;
\draw  [fill={rgb, 255:red, 0; green, 0; blue, 0 }  ,fill opacity=1 ] (142.16,159.54) .. controls (142.16,158.76) and (142.78,158.12) .. (143.53,158.12) .. controls (144.28,158.12) and (144.89,158.76) .. (144.89,159.54) .. controls (144.89,160.32) and (144.28,160.96) .. (143.53,160.96) .. controls (142.78,160.96) and (142.16,160.32) .. (142.16,159.54) -- cycle ;
\draw  [fill={rgb, 255:red, 0; green, 0; blue, 0 }  ,fill opacity=1 ] (161.84,164.24) .. controls (161.84,163.46) and (162.45,162.82) .. (163.21,162.82) .. controls (163.96,162.82) and (164.57,163.46) .. (164.57,164.24) .. controls (164.57,165.02) and (163.96,165.66) .. (163.21,165.66) .. controls (162.45,165.66) and (161.84,165.02) .. (161.84,164.24) -- cycle ;
\draw    (181.3,168.23) -- (181.3,176.44) ;
\draw    (181.3,168.23) -- (167.5,168.46) ;
\draw    (167.5,168.46) -- (167.5,176.44) ;
\draw    (163.21,164.24) -- (167.5,168.46) ;
\draw    (185.69,163.74) -- (181.3,168.23) ;
\draw    (137.64,153.71) -- (116.33,153.84) ;
\draw    (137.54,154.44) -- (116.23,154.57) ;
\draw    (137.54,154.44) -- (149.4,163.94) ;
\draw    (142.03,149.22) -- (137.64,153.71) ;
\draw    (149.4,163.94) -- (149.6,176.26) ;
\draw    (163.21,164.24) -- (149.4,163.94) ;
\draw    (163.17,159.14) -- (163.21,164.24) ;
\draw    (163.17,159.14) -- (116.93,159.48) ;
\draw    (167.57,154.65) -- (163.17,159.14) ;

\draw  [fill={rgb, 255:red, 0; green, 0; blue, 0 }  ,fill opacity=1 ] (122.12,204.14) .. controls (122.12,203.36) and (122.74,202.72) .. (123.49,202.72) .. controls (124.24,202.72) and (124.85,203.36) .. (124.85,204.14) .. controls (124.85,204.92) and (124.24,205.56) .. (123.49,205.56) .. controls (122.74,205.56) and (122.12,204.92) .. (122.12,204.14) -- cycle ;
\draw  [fill={rgb, 255:red, 0; green, 0; blue, 0 }  ,fill opacity=1 ] (179.58,218.23) .. controls (179.58,217.45) and (180.19,216.81) .. (180.94,216.81) .. controls (181.69,216.81) and (182.31,217.45) .. (182.31,218.23) .. controls (182.31,219.01) and (181.69,219.65) .. (180.94,219.65) .. controls (180.19,219.65) and (179.58,219.01) .. (179.58,218.23) -- cycle ;
\draw  [fill={rgb, 255:red, 0; green, 0; blue, 0 }  ,fill opacity=1 ] (141.8,209.54) .. controls (141.8,208.76) and (142.41,208.12) .. (143.17,208.12) .. controls (143.92,208.12) and (144.53,208.76) .. (144.53,209.54) .. controls (144.53,210.32) and (143.92,210.96) .. (143.17,210.96) .. controls (142.41,210.96) and (141.8,210.32) .. (141.8,209.54) -- cycle ;
\draw  [fill={rgb, 255:red, 0; green, 0; blue, 0 }  ,fill opacity=1 ] (161.48,214.24) .. controls (161.48,213.46) and (162.09,212.82) .. (162.85,212.82) .. controls (163.6,212.82) and (164.21,213.46) .. (164.21,214.24) .. controls (164.21,215.02) and (163.6,215.66) .. (162.85,215.66) .. controls (162.09,215.66) and (161.48,215.02) .. (161.48,214.24) -- cycle ;
\draw    (180.66,214.47) -- (180.36,226.54) ;
\draw    (181.47,215.13) -- (152.56,215.06) ;
\draw    (181.47,215.13) -- (181.27,226.46) ;
\draw    (185.05,209.97) -- (180.66,214.47) ;
\draw    (137.28,203.72) -- (115.97,203.85) ;
\draw    (137.18,204.45) -- (115.87,204.58) ;
\draw    (137.18,204.45) -- (143.17,209.54) ;
\draw    (141.67,199.22) -- (137.28,203.72) ;
\draw    (143.17,209.54) -- (115.93,209.44) ;
\draw    (180.66,214.47) -- (152.69,214.34) ;
\draw    (152.56,215.06) -- (152.72,226.19) ;
\draw    (152.69,214.34) -- (143.17,209.54) ;
\draw    (185.86,210.63) -- (181.47,215.13) ;

\draw  [fill={rgb, 255:red, 0; green, 0; blue, 0 }  ,fill opacity=1 ] (122.49,254.14) .. controls (122.49,253.36) and (123.1,252.72) .. (123.85,252.72) .. controls (124.6,252.72) and (125.22,253.36) .. (125.22,254.14) .. controls (125.22,254.92) and (124.6,255.56) .. (123.85,255.56) .. controls (123.1,255.56) and (122.49,254.92) .. (122.49,254.14) -- cycle ;
\draw  [fill={rgb, 255:red, 0; green, 0; blue, 0 }  ,fill opacity=1 ] (179.94,268.23) .. controls (179.94,267.45) and (180.55,266.81) .. (181.3,266.81) .. controls (182.06,266.81) and (182.67,267.45) .. (182.67,268.23) .. controls (182.67,269.01) and (182.06,269.65) .. (181.3,269.65) .. controls (180.55,269.65) and (179.94,269.01) .. (179.94,268.23) -- cycle ;
\draw  [fill={rgb, 255:red, 0; green, 0; blue, 0 }  ,fill opacity=1 ] (142.16,259.54) .. controls (142.16,258.76) and (142.78,258.12) .. (143.53,258.12) .. controls (144.28,258.12) and (144.89,258.76) .. (144.89,259.54) .. controls (144.89,260.32) and (144.28,260.96) .. (143.53,260.96) .. controls (142.78,260.96) and (142.16,260.32) .. (142.16,259.54) -- cycle ;
\draw  [fill={rgb, 255:red, 0; green, 0; blue, 0 }  ,fill opacity=1 ] (161.84,264.24) .. controls (161.84,263.46) and (162.45,262.82) .. (163.21,262.82) .. controls (163.96,262.82) and (164.57,263.46) .. (164.57,264.24) .. controls (164.57,265.02) and (163.96,265.66) .. (163.21,265.66) .. controls (162.45,265.66) and (161.84,265.02) .. (161.84,264.24) -- cycle ;
\draw    (181.3,268.23) -- (181.3,276.44) ;
\draw    (181.3,268.23) -- (167.5,268.46) ;
\draw    (167.5,268.46) -- (167.5,276.44) ;
\draw    (163.21,264.24) -- (167.5,268.46) ;
\draw    (185.69,263.74) -- (181.3,268.23) ;
\draw    (142.39,253.8) -- (133.01,253.91) -- (115.26,254.02) ;
\draw    (142.35,259.78) -- (115.54,259.86) ;
\draw    (142.35,259.78) -- (148.22,264.18) ;
\draw    (142.39,253.8) -- (142.35,259.78) ;
\draw    (148.22,264.18) -- (148.42,276.5) ;
\draw    (163.21,264.24) -- (148.22,264.18) ;
\draw    (163.11,254.37) -- (163.21,264.24) ;
\draw    (163.11,254.37) -- (116.24,254.74) ;
\draw    (166.21,251.14) -- (164.42,253) -- (163.11,254.37) ;
\draw    (146.49,250.58) -- (144.7,252.44) -- (142.39,253.8) ;

\draw  [fill={rgb, 255:red, 0; green, 0; blue, 0 }  ,fill opacity=1 ] (247.29,154.37) .. controls (247.29,153.59) and (247.9,152.96) .. (248.66,152.96) .. controls (249.41,152.96) and (250.02,153.59) .. (250.02,154.37) .. controls (250.02,155.16) and (249.41,155.79) .. (248.66,155.79) .. controls (247.9,155.79) and (247.29,155.16) .. (247.29,154.37) -- cycle ;
\draw  [fill={rgb, 255:red, 0; green, 0; blue, 0 }  ,fill opacity=1 ] (304.74,168.47) .. controls (304.74,167.68) and (305.36,167.05) .. (306.11,167.05) .. controls (306.86,167.05) and (307.47,167.68) .. (307.47,168.47) .. controls (307.47,169.25) and (306.86,169.88) .. (306.11,169.88) .. controls (305.36,169.88) and (304.74,169.25) .. (304.74,168.47) -- cycle ;
\draw  [fill={rgb, 255:red, 0; green, 0; blue, 0 }  ,fill opacity=1 ] (266.97,159.78) .. controls (266.97,158.99) and (267.58,158.36) .. (268.34,158.36) .. controls (269.09,158.36) and (269.7,158.99) .. (269.7,159.78) .. controls (269.7,160.56) and (269.09,161.19) .. (268.34,161.19) .. controls (267.58,161.19) and (266.97,160.56) .. (266.97,159.78) -- cycle ;
\draw  [fill={rgb, 255:red, 0; green, 0; blue, 0 }  ,fill opacity=1 ] (286.65,164.47) .. controls (286.65,163.69) and (287.26,163.06) .. (288.01,163.06) .. controls (288.77,163.06) and (289.38,163.69) .. (289.38,164.47) .. controls (289.38,165.26) and (288.77,165.89) .. (288.01,165.89) .. controls (287.26,165.89) and (286.65,165.26) .. (286.65,164.47) -- cycle ;
\draw    (306.11,168.47) -- (306.1,176.68) ;
\draw    (306.11,168.47) -- (292.3,168.69) ;
\draw    (292.3,168.69) -- (292.3,176.68) ;
\draw    (288.01,164.47) -- (292.3,168.69) ;
\draw    (310.5,163.97) -- (306.11,168.47) ;
\draw    (262.45,153.95) -- (241.14,154.08) ;
\draw    (262.35,154.68) -- (241.04,154.81) ;
\draw    (262.35,154.68) -- (268.34,159.78) ;
\draw    (266.84,149.45) -- (262.45,153.95) ;
\draw    (268.34,159.78) -- (268.33,176.35) ;
\draw    (288.01,164.47) -- (241.28,164.42) ;
\draw    (287.98,159.38) -- (288.01,164.47) ;
\draw    (287.98,159.38) -- (268.34,159.78) ;
\draw    (292.37,154.88) -- (287.98,159.38) ;

\draw  [fill={rgb, 255:red, 0; green, 0; blue, 0 }  ,fill opacity=1 ] (245.54,204.14) .. controls (245.54,203.36) and (246.15,202.72) .. (246.91,202.72) .. controls (247.66,202.72) and (248.27,203.36) .. (248.27,204.14) .. controls (248.27,204.92) and (247.66,205.56) .. (246.91,205.56) .. controls (246.15,205.56) and (245.54,204.92) .. (245.54,204.14) -- cycle ;
\draw  [fill={rgb, 255:red, 0; green, 0; blue, 0 }  ,fill opacity=1 ] (302.99,218.23) .. controls (302.99,217.45) and (303.6,216.81) .. (304.36,216.81) .. controls (305.11,216.81) and (305.72,217.45) .. (305.72,218.23) .. controls (305.72,219.01) and (305.11,219.65) .. (304.36,219.65) .. controls (303.6,219.65) and (302.99,219.01) .. (302.99,218.23) -- cycle ;
\draw  [fill={rgb, 255:red, 0; green, 0; blue, 0 }  ,fill opacity=1 ] (265.22,209.54) .. controls (265.22,208.76) and (265.83,208.12) .. (266.59,208.12) .. controls (267.34,208.12) and (267.95,208.76) .. (267.95,209.54) .. controls (267.95,210.32) and (267.34,210.96) .. (266.59,210.96) .. controls (265.83,210.96) and (265.22,210.32) .. (265.22,209.54) -- cycle ;
\draw  [fill={rgb, 255:red, 0; green, 0; blue, 0 }  ,fill opacity=1 ] (284.9,214.24) .. controls (284.9,213.46) and (285.51,212.82) .. (286.26,212.82) .. controls (287.02,212.82) and (287.63,213.46) .. (287.63,214.24) .. controls (287.63,215.02) and (287.02,215.66) .. (286.26,215.66) .. controls (285.51,215.66) and (284.9,215.02) .. (284.9,214.24) -- cycle ;
\draw    (304.36,218.23) -- (304.31,225.53) ;
\draw    (286.26,214.24) -- (276.11,214.33) ;
\draw    (286.26,214.24) -- (286.4,225.93) ;
\draw    (308.75,213.74) -- (304.36,218.23) ;
\draw    (260.7,203.71) -- (239.38,203.84) ;
\draw    (260.6,204.44) -- (239.28,204.57) ;
\draw    (260.6,204.44) -- (266.59,209.54) ;
\draw    (265.09,199.22) -- (260.7,203.71) ;
\draw    (266.59,209.54) -- (239.34,209.44) ;
\draw    (304.36,218.23) -- (280.27,218.22) ;
\draw    (276.11,214.33) -- (280.27,218.22) ;
\draw    (276.11,214.33) -- (266.59,209.54) ;
\draw    (290.65,209.74) -- (286.26,214.24) ;
\draw    (280.27,218.22) -- (280.27,225.46) ;

\draw  [fill={rgb, 255:red, 0; green, 0; blue, 0 }  ,fill opacity=1 ] (245.69,254.14) .. controls (245.69,253.36) and (246.3,252.72) .. (247.06,252.72) .. controls (247.81,252.72) and (248.42,253.36) .. (248.42,254.14) .. controls (248.42,254.92) and (247.81,255.56) .. (247.06,255.56) .. controls (246.3,255.56) and (245.69,254.92) .. (245.69,254.14) -- cycle ;
\draw  [fill={rgb, 255:red, 0; green, 0; blue, 0 }  ,fill opacity=1 ] (303.14,268.23) .. controls (303.14,267.45) and (303.75,266.81) .. (304.51,266.81) .. controls (305.26,266.81) and (305.87,267.45) .. (305.87,268.23) .. controls (305.87,269.01) and (305.26,269.65) .. (304.51,269.65) .. controls (303.75,269.65) and (303.14,269.01) .. (303.14,268.23) -- cycle ;
\draw  [fill={rgb, 255:red, 0; green, 0; blue, 0 }  ,fill opacity=1 ] (265.37,259.54) .. controls (265.37,258.76) and (265.98,258.12) .. (266.73,258.12) .. controls (267.49,258.12) and (268.1,258.76) .. (268.1,259.54) .. controls (268.1,260.32) and (267.49,260.96) .. (266.73,260.96) .. controls (265.98,260.96) and (265.37,260.32) .. (265.37,259.54) -- cycle ;
\draw  [fill={rgb, 255:red, 0; green, 0; blue, 0 }  ,fill opacity=1 ] (285.05,264.24) .. controls (285.05,263.46) and (285.66,262.82) .. (286.41,262.82) .. controls (287.17,262.82) and (287.78,263.46) .. (287.78,264.24) .. controls (287.78,265.02) and (287.17,265.66) .. (286.41,265.66) .. controls (285.66,265.66) and (285.05,265.02) .. (285.05,264.24) -- cycle ;
\draw    (304.51,268.23) -- (304.45,275.53) ;
\draw    (286.41,264.24) -- (276.26,264.33) ;
\draw    (286.41,264.24) -- (290.3,268.15) ;
\draw    (308.9,263.74) -- (304.51,268.23) ;
\draw    (260.85,253.71) -- (239.53,253.84) ;
\draw    (260.74,254.44) -- (239.43,254.57) ;
\draw    (260.74,254.44) -- (266.73,259.54) ;
\draw    (265.24,249.22) -- (260.85,253.71) ;
\draw    (266.73,259.54) -- (239.49,259.44) ;
\draw    (304.51,268.23) -- (290.3,268.15) ;
\draw    (276.26,264.33) -- (276.21,275.53) ;
\draw    (276.26,264.33) -- (266.73,259.54) ;
\draw    (290.8,259.74) -- (286.41,264.24) ;
\draw    (290.3,268.15) -- (290.3,275.4) ;

\node at (150,185) {$\displaystyle \langle 2\rangle +\langle 2d_{2} \rangle $};
\node at (275,185) {$\displaystyle \langle 1 \rangle $};

\node at (150,235) {$\displaystyle \langle 2d_{3} \rangle +\langle 2d_{4} \rangle $};
\node at (275,235) {$\displaystyle \langle 1 \rangle $};

\node at (150,285) {$\displaystyle \langle 2\rangle +\langle 2d_{1} \rangle $};
\node at (275,285) {$\displaystyle 2\mathbb{H} $};
\end{tikzpicture}

\end{center} \caption{The count of $N_{\AA^1,d}(\sigma)$ for $\sigma=\big(k(\sqrt{d_1}),k(\sqrt{d_2}),k(\sqrt{d_3}),k(\sqrt{d_4})\big)$ and $d=3$ using tropical stable maps. The point conditions are depicted as larger points, since they are double. The picture indicates the parametrization by the abstract graphs $\Gamma$ for each of the tropical stable maps $(\Gamma,f)$ as well as the image $f(\Gamma)$ in such a way that edges that actually map on top of each other are drawn closeby. These pictures can be obtained from the tropical stable maps in Figure \ref{fig-ninecubics} by merging pairs of points.}
    \label{fig-tropQEQ}
\end{figure}
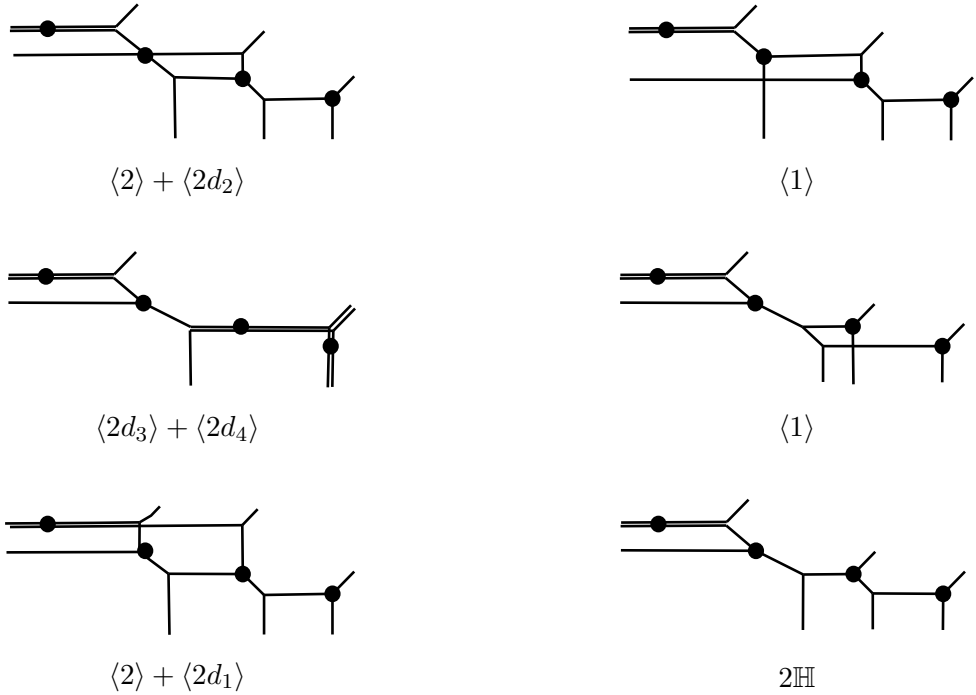

\begin{example}\label{ex-tropQEQ}
    Let $d=3$, and $\sigma=\big(k(\sqrt{d_1}),k(\sqrt{d_2}),k(\sqrt{d_3}),k(\sqrt{d_4})\big)$. The count $N_{\AA^1,d}(\sigma)$ can be determined using a quadratically enriched count of tropical stable maps satisfying double point conditions, as depicted in Figure \ref{fig-tropQEQ}.
    Using $2\gw{2} = 2\gw{1}$ in $\GW(k)$ we obtain $N_{\AA^1,3}(\sigma)=4\cdot \langle 1\rangle +\sum_{i=1}^4\langle 2d_{i} \rangle+ 2\mathbb{H}$.

Notice that this invariant encompasses as many as six real and complex counts of curves (in addition to various other quadratically enriched counts):
\begin{itemize}
\item If we specialize to the ground field $k=\mathbb{C}$, i.e.\ take the rank, we obtain the Gromov-Witten invariant $N_3=12$.
\item If we specialize to $k=\mathbb{R}$ and take the signature, we obtain the Welschinger invariants: 
\begin{itemize}
    \item Setting $d_i=1$ for all  $i$ (i.e.\ we do not really work in a quadratic extension but just in $\mathbb{R}$), we obtain $W_3=8$.
    \item For $l = 1, ..., 4$, setting the first $4-l$ of the $d_i = 1$ and the remaining $d_i = -1$, then we obtain $W_{3, l} = 8 - 2l$.
\end{itemize}
\item If we merely set some of the $d_i$ to 1, then we specialize the counts for other $\sigma$'s, e.g. set $d_3 = d_4 = 1$ and get the count for $\sigma' = (k(\sqrt{d_1}), k(\sqrt{d_2}), k, k, k, k)$, which is $N_{\AA^1, 3}(\sigma') = 6 \gw{1} + \gw{2d_1} + \gw{2d_2} + 2 \mathbb{H}$.
\end{itemize}    
\end{example}

\begin{remark}
While the general idea of a correspondence theorem enabling us to compute $N_{\AA^1,d}(\sigma)$ by means of tropical geometry is the same as described in Sections \ref{sec-corres}, \ref{sec-tropW} and \ref{sec-tropQE}, the details get more and more involved. On the one hand, the underlying combinatorics are way more involved by the need to merge points into double points. On the other hand, the construction of local pieces of curves and the determination of their nodes and quadratically enriched Welschinger signs is also way more complicated. In \cite{JPMPR25}, we have therefore restricted ourselves to special point configurations in the tropical world and this is reflected by our choice of point configuration in Example \ref{ex-tropQEQ}. Since we know that the count is invariant, this does not limit the power of our result to compute quadratically enriched invariants. Nonetheless, it would be interesting to establish a full correspondence theorem for the count involving quadratic field extensions.
\end{remark}

\begin{remark}\label{rem-otherfieldextensions}
    If we want to allow point conditions defined over field extensions of higher degree, we would have to merge even more points in the tropical world. We expect that both the underlying combinatorics as well as the constructions for local pieces in a possible correspondence are even more difficult, but we still see it as an interesting task for further research to establish such tropical counts.
\end{remark}

\bibliographystyle{plain} 
\bibliography{bibliographie} 
\end{document}